\documentclass[12pt]{amsart}

\usepackage{amssymb}
\usepackage{amscd}
\usepackage[all]{xy}
\swapnumbers
\def\frak{\mathfrak}
\def\Bbb{\mathbb}
\def\Cal{\mathcal}

\newtheorem*{prop*}{Proposition}

\newtheorem*{thm*}{Theorem}

\newtheorem*{lem*}{Lemma}

\newtheorem*{kor*}{Corollary}

\newcommand{\id}{\operatorname{id}}

\renewcommand{\ker}{\operatorname{ker}}
\newcommand{\im}{\operatorname{im}}

\newcommand{\x}{\times}
\renewcommand{\o}{\circ}

\let\ccdot\cdot
\def\cdot{\hbox to 2.5pt{\hss$\ccdot$\hss}}

\newcommand{\g}{{\frak g}}

\def\gog{\frak{g}}

\newcommand{\calH}{{\Cal H}}

\newcommand{\bV}{{\Bbb V}}

\newcommand{\bE}{{\Bbb E}}
\newcommand{\bF}{{\Bbb F}}

\newcommand{\bC}{{\Bbb C}}

\newcommand{\al}{\alpha}
\newcommand{\be}{\beta}
\newcommand{\ga}{\gamma}
\newcommand{\de}{\delta}

\newcommand{\ka}{\kappa}
\newcommand{\la}{\lambda}
\newcommand{\om}{\omega}
\renewcommand{\phi}{\varphi}
\newcommand{\ph}{\varphi}
\newcommand{\ps}{\psi}
\newcommand{\si}{\sigma}

\newcommand{\De}{\Delta}
\newcommand{\Ga}{\Gamma}
\newcommand{\La}{\Lambda}
\newcommand{\Om}{\Omega}
\newcommand{\Ph}{\Phi}

\newcommand{\Si}{\Sigma}

\def\sideremark#1{\ifvmode\leavevmode\fi\vadjust{%            The remark
\vbox to0pt{\hbox to 0pt{\hskip\hsize\hskip1em%               will appear only
\vbox{\hsize3cm\tiny\raggedright\pretolerance10000%          on the side
\noindent #1\hfill}\hss}\vbox to8pt{\vfil}\vss}}}%           in 3cm  
\newcommand\bxbdb[4]{\begin{picture}(96,12)\put(5,3){\line(1,0){36}}%
\put(45,3){\line(1,0){6}}\put(79,3){\line(-1,0){6}}%
\put(63,3){\makebox(0,0){\dots}}%
\put(3,3){\makebox(0,0){$\o$}}\put(23,3){\makebox(0,0){$\x$}}%
\put(43,3){\makebox(0,0){$\o$}}\put(81,3){\makebox(0,0){$\o$}}%
\put(3,10){\makebox(0,0){\scriptsize $#1$}}%
\put(23,10){\makebox(0,0){\scriptsize $#2$}}%
\put(43,10){\makebox(0,0){\scriptsize $#3$}}%
\put(81,10){\makebox(0,0){\scriptsize $#4$}}\end{picture}}

\newcommand\xbdbx[5]{\begin{picture}(76,12)\put(3,3){\line(1,0){18}}%
\put(25,3){\line(1,0){6}}\put(73,3){\line(-1,0){18}}%
\put(51,3){\line(-1,0){6}}\put(39,3){\makebox(0,0){\dots}}%
\put(3,3){\makebox(0,0){$\x$}}\put(23,3){\makebox(0,0){$\o$}}%
\put(53,3){\makebox(0,0){$\o$}}\put(73,3){\makebox(0,0){$\x$}}%
\put(3,10){\makebox(0,0){\scriptsize $#1$}}%
\put(23,10){\makebox(0,0){\scriptsize $#2$}}%
\put(38,10){\makebox(0,0){\scriptsize $#3$}}%
\put(53,10){\makebox(0,0){\scriptsize $#4$}}%
\put(73,10){\makebox(0,0){\scriptsize $#5$}}\end{picture}}

\newcommand\Cbxbdbb[5]{\begin{picture}(96,12)\put(5,3){\line(1,0){36}}%
\put(45,3){\line(1,0){6}}%
\put(71,3){\line(-1,0){6}}\put(74,4){\line(1,0){18}}%
\put(74,2){\line(1,0){18}}\put(59,3){\makebox(0,0){\dots}}%
\put(83,3){\makebox(0,0){$<$}}%
\put(3,3){\makebox(0,0){$\o$}}\put(23,3){\makebox(0,0){$\x$}}\put(43,3){\makebox(0,0){$\o$}}
\put(73,3){\makebox(0,0){$\o$}}\put(93,3){\makebox(0,0){$\o$}}%
\put(3,10){\makebox(0,0){\scriptsize $#1$}}%
\put(23,10){\makebox(0,0){\scriptsize $#2$}}%
\put(43,10){\makebox(0,0){\scriptsize $#3$}}%
\put(73,10){\makebox(0,0){\scriptsize $#4$}}%
\put(93,10){\makebox(0,0){\scriptsize $#5$}}\end{picture}}

\begin{document}
\title{Subcomplexes in Curved BGG--Sequences} 
\author{Andreas \v Cap and Vladim\'{i}r Sou\v cek}

\address{A.C.: Institut f\"ur Mathematik, Universit\"at Wien, Strudlhofgasse 4,
A--1090 Wien, Austria and International Erwin Schr\"odinger Institute for
Mathematical Physics, Boltzmanngasse 9, A--1090 Wien, Austria\newline\indent
V.S.: Mathematical Institute, Charles University, Sokolovsk\'a 83,
Praha, Czech Republic}
\email{Andreas.Cap@esi.ac.at, soucek@karlin.mff.cuni.cz}
\begin{abstract}
BGG--sequences offer a uniform construction for invariant differential
operators for a large class of geometric structures called parabolic
geometries. For locally flat geometries, the resulting sequences are
complexes, but in general the compositions of the operators in such a
sequence are nonzero. In this paper, we show that under appropriate
torsion freeness and/or semi--flatness assumptions certain parts of
all BGG sequences are complexes.

Several examples of structures, including quaternionic structures,
hypersurface type CR structures and quaternionic contact structures
are discussed in detail. In the case of quaternionic structures we
show that several families of complexes obtained in this way are
elliptic.  
\end{abstract}

\maketitle

\section{Introduction}\label{1}
Parabolic geometries form a large class of geometric structures
containing examples like conformal, quaternionic, hypersurface type
CR, and certain higher codimension CR structures. Via the
interpretation as Cartan geometries with homogeneous model a
generalized flag manifold, these structures can be studied in a
surprisingly uniform way. An important and difficult problem is the
construction of invariant differential operators for such geometries,
i.e.~operators which are intrinsic to the structure. 

For the homogeneous model $G/P$ (and geometries locally isomorphic to
the homogeneous model) this problem can be reformulated in terms of
representation theory. Via the theory of homomorphisms of generalized
Verma modules one obtains an almost complete answer. In particular,
invariant differential operators between sections of bundles
associated to irreducible representations show up in patterns which can
be described combinatorially in terms of the Weyl group of the Lie
algebra $\frak g$ of $G$. The different patterns are indexed by
certain weights for $\frak g$. For dominant integral weights (which
covers most of the cases of interest), the resulting pattern forms a
resolution of the finite dimensional irreducible representation of
$\frak g$, the celebrated (generalized) Bernstein--Gelfand--Gelfand (or
BGG) resolution, see \cite{BGG, Lepowsky}.

The representation theory arguments leading to the BGG resolutions are
combinatorial in nature and finding an explicit interpretation in
terms of differential operators is a highly nontrivial task. An
independent construction of the BGG resolutions in terms of
differential operators was given in \cite{CSSBGG} and improved in
\cite{CD}. This construction has the advantage that it works without
changes for curved geometries, thus providing a construction of a
large number of invariant differential operators for arbitrary
parabolic geometries. The resulting patterns of operators are referred
to as BGG--sequences. The construction also relates sections of the
bundles showing up in a BGG--sequence to differential forms with
values in a so--called tractor bundle and the (higher order) operators
in the sequence to a covariant exterior derivative on these forms.

A drawback of the curved BGG--sequences is that the operators in the
sequence have nontrivial compositions in general, so usually one does
not obtain complexes in this way. The machinery of \cite{CSSBGG, CD}
however is strong enough to give explicit formulae for the
compositions. Starting from these formulae, we prove a simple
criterion (in terms of weights), which ensures that some of the
compositions do vanish provided that the harmonic curvature of the
geometry satisfies certain restrictions, see Theorem \ref{3.3}. The
necessary restrictions on the curvature usually include torsion
freeness, but in some cases one also needs assumptions like
semi--flatness, the most prominent of those being (anti)self duality
in four dimensional conformal geometry.

Using finite dimensional representation theory, the weight conditions
are then systematically studied in several examples. We describe the
form of the BGG patterns and identify in each case several
sub--patterns, for which the weight conditions are always satisfied.
Under the appropriate curvature restrictions, these sub--patterns
therefore give rise to subcomplexes in each BGG sequence.

The first examples we discuss are almost Grassmannian and almost
quaternionic structures, which are different real forms of the same
complex geometry. In dimension four, the structures can be
equivalently described as conformal structures in split signature
respectively Riemannian signature. The curvature conditions amount to
torsion freeness for higher dimensions and (anti) self duality in
dimension four. Our results in this case vastly generalize the
complexes found in \cite{Salamon} and \cite{Baston} for
quaternionic structures. We should also mention that many known
applications of curved BGG--sequences actually use known special cases
of these subcomplexes, see \cite{David}.

Next, we discuss Lagrangean contact structures and partially
integrable almost CR structures of hypersurface type, which again are
different real forms of the same complex geometry. The appropriate
curvature restriction again is torsion freeness, which is equivalent
to integrability in the CR case. 

Finally, we discuss the case of quaternionic contact structures as
introduced by O.~Biquard, see \cite{Biq,Biq2}. The curvature condition
is torsion freeness in the case of dimension seven, while it is
automatically satisfied in higher dimensions. 

The complexes for quaternionic structures found in \cite{Salamon}
(which involve first order operators only) are elliptic. In the last
section we extend this result to other families of subcomplexes in the
quaternionic case, which involve operators of arbitrarily high orders.

The BGG sequence associated to the adjoint representation is closely
related to the theory of infinitesimal deformations of parabolic
geometries. For the examples of structures discussed in this paper,
one of the subcomplexes in the adjoint BGG sequence can be naturally
interpreted as a deformation complex in the subcategory of structures
satisfying the curvature restrictions. This is discussed in
\cite{deformations}.

\noindent{\bf Acknowledgments.}  The research evolved during the
meetings of authors at the Erwin Schr\"odinger Institute for
Mathematical Physics in Vienna and the Charles University in
Prague. First author supported by project P15747--N05 of the Fonds zur
F\"orderung der wissenschaftlichen Forschung (FWF). The second author
thanks the grant GA\v CR Nr.~201/02/1390 and the institutional grant
MSM 21620839 for their support.

\section{Background}\label{2}
In this section, we briefly review some basic facts about parabolic
geometries and BGG sequences, mainly to fix the notation used in the
sequel. More detailed information can be found in
\cite{C-tw,Weyl,CSSBGG}. 

\subsection{Parabolic geometries}\label{2.1}
A type of parabolic geometries is determined by a parabolic subgroup
$P$ in a semisimple Lie group $G$. Parabolic subgroups can be nicely
described in terms of $|k|$--gradings of the Lie algebra $\frak g$ of
$G$. Details about $|k|$--gradings can be found in \cite{Yam}. 

A $|k|$--grading on a semisimple Lie algebra $\frak g$ is a decomposition 
$$ 
\frak g=\frak g_{-k}\oplus\dots \oplus \frak g_{-1}\oplus\frak
g_0\oplus\frak g_1\oplus\dots\oplus \frak g_k
$$ such that $[\frak g_i,\frak g_j]\subset\frak g_{i+j}$ and such that
the subalgebra $\frak g_-:=\frak g_{-k}\oplus\dots\oplus\frak g_{-1}$
is generated by $\frak g_{-1}$. Defining $\frak g^i:=\frak
g_i\oplus\dots\oplus\frak g_k$, we obtain a filtration of $\frak g$
and $[\frak g^i,\frak g^j]\subset\frak g^{i+j}$. In particular, $\frak
p:=\frak g^0$ is a Lie subalgebra of $\frak g$ and $\frak p_+:=\frak
g^1$ is a nilpotent ideal in $\frak p$.

On the group level, we define $G_0\subset P\subset G$ as the subgroups
of elements whose adjoint action preserves each grading component
$\frak g_i$ respectively each filtration component $\frak g^i$. It
turns out that $\frak p\subset\frak g$ is a parabolic subalgebra and
$P=N_G(\frak p)$ is the usual parabolic subgroup associated to $\frak
p$, while $G_0$ has Lie algebra $\frak g_0$. Moreover, exponential map
defines a diffeomorphism from $\frak p_+$ onto a closed normal
subgroup $P_+\subset P$, and $P$ is the semidirect product of $G_0$ and
$P_+$. In the complex case, parabolic subalgebras (up to conjugacy)
are in bijective correspondence with sets of simple roots. This leads
to a description of $|k|$--gradings in terms of Dynkin diagrams with
crosses, see \cite{BE}. In the real case, there is a similar
description in terms of Satake diagrams, see \cite{Yam}.

Parabolic geometries are then defined as Cartan geometries of type
$(G,P)$. This means that a parabolic geometry of type $(G,P)$ on $M$
consists of a principal $P$--bundle $p:\Cal G\to M$ endowed with a
Cartan connection $\om\in\Om^1(\Cal G,\frak g)$. The \textit{homogeneous
  model} of the geometry is given by the natural bundle $p:G\to G/P$
endowed with the left Maurer-Cartan form as a Cartan connection. A
\textit{morphism} of parabolic geometries is a principal bundle map
which is compatible with the Cartan connections. In particular, any
morphism is a local diffeomorphism. 

The curvature of a Cartan connection $\om$ is defined as the $\frak
g$--valued two--form $K\in\Om^2(\Cal G,\frak g)$ defined by the
structure equation
$$
K(\xi,\eta)=d\om (\xi,\eta)+[\om(\xi),\om(\eta)],
$$ where $\xi$ and $\eta$ are vector fields on $\Cal G$ and the
bracket is in $\frak g$. The form $K$ is horizontal and equivariant,
so it may be interpreted as a two--form $\ka$ on $M$ with values in
the associated bundle $\Cal AM:=\Cal G\x_P\frak g$, the
\textit{adjoint tractor bundle}. The $P$--invariant filtration
$\{\frak g^i\}$ of $\frak g$ gives rise to a filtration $\Cal AM=\Cal
A^{-k}M\supset\dots\supset\Cal A^kM$ by smooth subbundles and the Lie
bracket on $\frak g$ gives rise to an algebraic bracket $\{\ ,\ \}$ on
$\Cal AM$ making it into a bundle of filtered Lie algebras modeled on
$\frak g$.

The Cartan connection $\om$ induces an isomorphism $TM\cong\Cal
G\x_P(\frak g/\frak p)$. Hence there is a natural projection $\Pi:\Cal
AM\to\Cal TM$ which induces an isomorphism $\Cal AM/\Cal A^0M\cong
TM$. Via this isomorphism, the filtration of $\Cal AM$ descends to a
filtration $TM=T^{-k}M\supset\dots\supset T^{-1}M$ of the tangent
bundle by smooth subbundles. Applying the projection $\Pi$ to the
values of $\ka$ we obtain a $TM$--valued two--form $\ka_-$, which is
called the \textit{torsion} of the Cartan connection $\om$. The
geometry is called \textit{torsion free} if this torsion vanishes.

Via the filtrations of $TM$ and $\Cal AM$ one has a natural notion of
homogeneity for $\Cal AM$--valued differential forms. In particular,
we say that $\ka$ is homogeneous of degree $\geq\ell$ if
$\ka(T^iM,T^jM)\subset\Cal A^{i+j+\ell}M$ for all $i,j=-k,\dots,-1$. A
parabolic geometry is called \textit{regular} if its curvature is
homogeneous of degree $\geq 1$. Note that torsion free parabolic
geometries are automatically regular.

\subsection{Lie algebra homology and normalization}\label{2.4}
Parabolic geometries are mainly interesting as an equivalent
conceptual description for a large number of (seemingly very diverse)
examples of geometric structures. Usually, the given geometric
structure can be easily encoded into what is called a regular
infinitesimal flag structure, see \cite{Weyl}. This consists of a
filtration $\{T^iM\}$ of $TM$ and a principal $G_0$--bundle $\Cal
G_0\to M$ endowed with certain partially defined differential forms.
Under a cohomological condition, which is satisfied for all the
structures considered in this paper, one can then apply involved
prolongation procedures (see \cite{Tan, Mor, CS}). These extend
$\Cal G_0$ to a principal $P$--bundle $p:\Cal G\to M$ endowed with a
Cartan connection $\om$. In particular, the given filtration of $TM$
coincides with the one obtained from $(\Cal G,\om)$ as in \ref{2.1}
and $\Cal G/P_+\cong\Cal G_0$.  The resulting parabolic geometry is
uniquely determined (up to isomorphism) if one in addition requires
the curvature of $\om$ to satisfy a normalization condition to be
discussed below. This leads to an equivalence of categories between
regular normal parabolic geometries and the underlying structures.

By forming associated bundles to the Cartan bundle, natural vector
bundles on parabolic geometries of type $(G,P)$ are defined by
representations of the Lie group $P$. A simple class of such
representations are obtained by trivially extending representations of
$G_0$. In particular, all irreducible representations of $P$ are
obtained in this way. Since the group $G_0$ turns out to be reductive,
its representation theory is well understood. If $\Bbb E$ is
such a representation, then $\Cal G\x_P\Bbb E\cong\Cal G_0\x_{G_0}\Bbb
E$. Therefore, the corresponding natural vector bundles admit a direct
interpretation in terms of the underlying geometric structure.

A second source of representations of $P$ is restrictions of
representations of $G$. The corresponding natural vector bundles are
called tractor bundles. While from a geometrical point of view these
are rather unusual objects, they are in some respects easy to deal
with. For example, the Cartan connection $\om$ induces a linear
connection on each tractor bundle. The general theory of tractor
bundles is developed in \cite{Cap-Gover}.

An important link between the two classes of representations we have
discussed is given by Lie algebra homology. Let $\Bbb V$ be a
representation of $G$, viewed as a representation of $P$ by
restriction. Infinitesimally, we in particular get a representation of
$\frak p_+$ on $\Bbb V$. The standard complex for computing the Lie
algebra homology $H_*(\frak p_+,\Bbb V)$ has the form
$$
0\to \Bbb V\overset{\partial^*}{\longrightarrow} \frak p_+\otimes\Bbb
V\to \dots \to \La^k\frak p_+\otimes\Bbb
V\overset{\partial^*}{\longrightarrow} \La^{k+1}\frak p_+\otimes\Bbb
V\to\dots 
$$
Following the literature on parabolic geometries, we denote the
standard differential by $\partial^*$ and call it the \textit{Kostant
  codifferential}. Explicitly, it is given by
\begin{align*}
  \partial^*(Z_0\wedge&\dots\wedge Z_n\otimes
  v)=\textstyle\sum_{i=0}^n(-1)^{i+1}
  Z_0\wedge\cdots\hat i\dots\wedge Z_n\otimes Z_i\cdot v+\\
  &+\textstyle\sum_{i<j}(-1)^{i+j}[Z_i,Z_j]\wedge\cdots\hat
  i\cdots\hat j\dots\wedge Z_n\otimes v
\end{align*}
with hats denoting omission.

Evidently, all spaces in the standard complex are representations of
$P$ and $\partial^*$ is $P$--equivariant. Hence each of the homology
groups $H_k(\frak p_+,\Bbb V)$ naturally is a $P$--module. Moreover,
one easily verifies that $P_+$ acts trivially on the homology, so the
representations $H_k(\frak p_+,\Bbb V)$ come from the subgroup $G_0$.

In \cite{Kostant}, B.~Kostant gave an explicit algorithm to compute,
for each $k$, the $G_0$--module $H^k(\frak p_+,\Bbb V)$, which is dual
to $H_k(\frak p_+,\Bbb V^*)$, in the case when $\frak g$ is complex
and simple and $\Bbb V$ is a complex irreducible representation. Using
basic tricks of the trade in Lie algebra cohomology one may also deal
with the real cases, so all the modules in question are explicitly
computable.

This construction has a direct geometric counterpart. We have noted
above that $TM\cong\Cal G\x_P\frak g/\frak p$, with the action coming
from the adjoint representation. Now $(\frak g/\frak p)^*\cong\frak
p_+$ as $P$--modules via the Killing form of $\frak g$. Hence
$T^*M\cong\Cal G\x_P\frak p_+=\Cal A^1M$, so $T^*M$ naturally is a
bundle of nilpotent filtered Lie algebras with the restriction of the
bracket $\{\ ,\ \}$ from \ref{2.1}. Denoting the tractor bundle
corresponding to the representation $\Bbb V$ by $VM$, we see that
$\Cal G\x_P(\La^k\frak p_+\otimes\Bbb V)\cong \La^kT^*M\otimes VM$,
and the codifferential induces a natural bundle map
$$
\partial^*:\La^kT^*M\otimes VM\to\La^{k-1}T^*M\otimes VM.
$$
Of course, we have $\partial^*\o\partial^*=0$ and the kernel and image
of $\partial^*$ are natural subbundles. Their quotient is by
construction isomorphic to $\Cal G\x_PH_k(\frak p_+,\Bbb V)$. We
denote this bundle by $H_k(T^*M,VM)$ since it is obtained from taking
the pointwise Lie algebra homology of $T_x^*M$ with coefficients in
the module $V_xM$. In particular, there is a natural bundle map
$\pi_H:\ker(\partial^*)\to H_k(T^*M,VM)$ and we denote by the same
symbol the induced tensorial map on sections.

The bundle map $\partial^*$ also induces a tensorial operator on
$VM$--valued differential forms, which we denote by the same symbol.
In the special case $\Bbb V=\frak g$ and $k=2$, we obtain
$\partial^*:\Om^2(M,\Cal AM)\to\Om^1(M,\Cal AM)$. A parabolic geometry
is called \textit{normal} if and only if $\partial^*(\ka)=0$, where
$\ka$ denotes the Cartan curvature. This is the normalization
condition referred to above. If this is satisfied then one defines the
\textit{harmonic curvature} $\ka_H:=\pi_H(\ka)\in\Ga(H_2(T^*M,\Cal
AM))$. This is much easier to interpret in terms of the underlying
structure than the full curvature $\ka$. In the regular normal case,
$\ka_H$ is a complete obstruction to local flatness.

\subsection{Strongly invariant operators}\label{2.3}
These form a class of invariant differential operators, which
algebraically are of particularly simple nature. The starting point
for this is that the first jet prolongation of a vector bundle
associated to $\Cal G$ can be identified with an associated bundle to
$\Cal G$: For an arbitrary representation $\Bbb E$ of $P$, one defines
$J^1\Bbb E:=\Bbb E\oplus L(\frak g/\frak p,\Bbb E)$ endowed with a
certain $P$--module structure. For each parabolic geometry $(p:\Cal
G\to M,\om)$, one has the natural vector bundle $EM:=\Cal G\x_P\Bbb
E$. One shows that the first jet prolongation $J^1EM$ of this bundle
is naturally isomorphic to $\Cal G\x_PJ^1\Bbb E$ for that $P$--module
structure. 

This does not extend to higher jets, but it does work for higher
semi--holonomic jet prolongations. One can put a $P$--module structure
on the space $\bar J^r\Bbb E=\oplus_{i=0}^r\otimes^i(\frak g/\frak
p)^*\otimes\Bbb E$, such that the corresponding associated bundle is
naturally isomorphic to the $r$th semi--holonomic jet prolongation
$\bar J^rEM$. The upshot of this is that any $P$--module homomorphisms
$\Psi:\bar J^r\Bbb E\to \Bbb F$ gives rise to a vector bundle map
$\bar J^rEM\to FM$, and thus to a natural $r$th order differential
operator $\Ga(EM)\to\Ga(FM)$. Operators arising in this way are called
strongly invariant. See \cite{CSSBGG, Sl97} for more information on
these issues.

\subsection{BGG sequences}\label{3.1}
We next sketch the geometric construction of the generalized BGG
resolutions as introduced in \cite{CSSBGG} and improved in \cite{CD}.
A nice explanation of the role of BGG sequences and some applications
can be found in \cite{MikeBGG,MikeNotices,David}. 

Let $\Bbb V$ be a representation of $G$ and for a parabolic geometry
$(p:\Cal G\to M,\om)$ let $VM$ be the corresponding tractor bundle. As
we have noted in \ref{2.4}, the Cartan connection $\om$ induces a
linear connection $\nabla=\nabla^{\Bbb V}$, called the tractor
connection, on $VM$. This extends to an operation
$$
d^\nabla:\Om^k(M,VM)\to\Om^{k+1}(M,VM)
$$
on $VM$--valued differential forms, called the covariant exterior
derivative. This operation is defined by taking the usual formula for
the exterior derivative and replacing the action of vector fields on
functions by the covariant derivative. 

The tractor bundle $VM$ carries a natural filtration by smooth
subbundles (see \cite{CSSBGG}) and correspondingly one has the notion
of homogeneity for $VM$--valued differential forms. The codifferential
$\partial^*$ from \ref{2.4} is compatible with homogeneities. For
regular normal parabolic geometries, also $d^\nabla$ is compatible
with homogeneities. Now one can view the composition $\partial^*\o
d^\nabla$ as an operator acting on sections of the bundle
$\im(\partial^*)\subset\La^kT^*M\otimes VM$. This preserves
homogeneities and one verifies that its homogeneous component of
degree zero is tensorial and invertible. (For this to have geometric
meaning one has to view the homogeneous component of degree zero as
acting on sections of the associated graded bundles.) This implies
that $\partial^*\o d^\nabla$ itself is invertible, and the inverse is
a (by construction natural) differential operator $Q$ acting on
sections of $\im(\partial^*)$. 

Now we define an operator $L:\Ga(H_k(T^*M,VM))\to\Om^k(M,VM)$ as
follows: For a section $\al$ of $H_k(T^*M,VM)$ choose a representative
$\ph\in\Om^k(M,VM)$, i.e.~$\partial^*(\ph)=0$ and $\pi_H(\ph)=\al$,
and put 
$$
L(\al):=\ph-Q\partial^*d^\nabla\ph.
$$ 
Since different choices for $\ph$ differ by sections of the subbundle
$\im(\partial^*)$ and the operator $Q\partial^*d^\nabla$ is the
identity on such sections, this is a well defined invariant operator,
called the \textit{splitting operator}. Let us collect its main
properties:

\begin{thm*}
  For any $\al\in\Ga(H_k(T^*M,VM))$ we have $\partial^*(L(\al))=0$,
  $\pi_H(L(\al))=\al$, and $\partial^*(d^\nabla L(\al))=0$. These
  three properties characterize $L(\al)$.
\end{thm*}
\begin{proof}
Choosing a representative $\ph$ for $\al$, we have
$L(\al)=\ph-Q\partial^*d^\nabla\ph$. Since $\partial^*(\ph)=0$ and $Q$
has values in $\im(\partial^*)$ we see that
$\partial^*(L(\al))=0$ and $\pi_H(L(\al))=\pi_H(\ph)=\al$. Since $Q$
is inverse to $\partial^*d^\nabla$ on $\im(\partial^*)$ we see that
$\partial^*d^\nabla L(\al)=\partial^*d^\nabla\ph-\partial^*d^\nabla\ph=0$. 

Conversely, assume that $\ps\in\Om^k(M,VM)$ satisfies
$\partial^*\ps=0$, $\pi_H(\ps)=\al$, and $\partial^*d^\nabla\ps=0$.
Then we can use $\ps$ as a representative for $\al$ in the
construction of $L$, and since $\partial^*d^\nabla\ps=0$ we get
$L(\al)=\ps$.
\end{proof}

Since $\partial^*d^\nabla L(\al)=0$, we obtain an invariant
differential operator
$$
D=D^\Bbb V:=\pi_H\o d^\nabla\o
L:\Ga(H_k(T^*M,VM))\to\Ga(H_{k+1}(T^*M,VM)),
$$
and these operators form the BGG sequence. 

It is well known that $d^\nabla\o d^\nabla$ is given by the action of
the curvature of $\nabla$. The curvature of a tractor connection is
given by the action of the Cartan curvature $\ka\in\Om^2(M,\Cal AM)$,
so we obtain $d^\nabla d^\nabla \ph=\ka\wedge\ph$. This is the
alternation of $(\xi_0,\dots,\xi_{k+1})\mapsto
\ka(\xi_0,\xi_1)\bullet\ph(\xi_2,\dots,\xi_{k+1})$ with the bundle map
$\bullet:\Cal AM\x VM\to VM$ induced by the action of $\frak g$ on
$\Bbb V$. For locally flat geometries, we have $\ka=0$ and the twisted
de-Rham sequence is a complex. This easily implies that also the BGG
sequence is a complex and both complexes compute the same cohomology,
see \cite{CSSBGG}. In the curved case, the compositions of the BGG
operators are nontrivial in general.

\section{Subcomplexes}\label{3}

\subsection{Compositions of BGG operators}\label{3.3}
Given a representation $\Bbb V$ of $G$, the representations $H_k(\frak
p_+,\Bbb V)$ of $G_0$ are always completely reducible. Hence they
split into a direct sum of irreducible components and correspondingly
the bundles $H_k(T^*M,VM)$ decompose into a direct sum of smooth
subbundles. Doing this for $k$ and $k+1$, the BGG operator splits into
a family of operators acting between sections of the individual
components. Likewise, the composition of two consecutive BGG operators
splits into components acting between the irreducible pieces. Assuming
restrictions on the Cartan curvature we can derive a purely algebraic
criterion for vanishing of pieces of the composition.

\begin{thm*}
  Let $\bE_0$ be a $G_0$--submodule of $H_2(\frak p_+,\frak g)$ and
  let $\bF_1$ and $\bF_2$ be irreducible components of $H_k(\frak
  p_+,\bV)$ and $H_{k+2}(\frak p_+,\bV)$, respectively. Suppose
  further that $\bF_2$ is not isomorphic to a $G_0$-submodule of
  $\oplus_{i\geq 0}\left(\otimes^i\frak p_+\otimes\bE_0\otimes
    \bF_1\right)$.
  
  Then for any torsion free normal parabolic geometry $(p:\Cal G\to
  M,\om)$ whose harmonic curvature $\ka_H$ is a section of
  $E_0M\subset H_2(T^*M,\Cal AM)$, the component in $F_2M\subset
  H_{k+2}(T^*M,VM)$ of the restriction of $D\o D$ to $F_1M\subset
  H_k(T^*M,VM)$ vanishes identically.
\end{thm*}
\begin{proof}
By torsion freeness, the covariant exterior derivative $d^\nabla$
coincides with the twisted exterior derivative $d^\Bbb V$ used in
\cite{CSSBGG}. Hence the constructions of the splitting operators and
BGG operators described in \ref{3.1} coincides with the construction
in \cite{CSSBGG}, so in particular all the operators are strongly
invariant. 

The Bianchi identity for linear connection implies $d^\nabla\ka=0$,
which together with $\partial^*(\ka)=0$ shows that $\ka=L(\ka_H)$. By
\cite[Theorem 2.5]{CSSBGG} the fact that $\ka_H\in\Ga(E_0M)$ implies
that $L(\ka_H)$ is a section $\Cal G\x_P\Bbb
E\subset\La^2T^*M\otimes\Cal AM$, where $\Bbb E\subset\La^2\frak
p_+\otimes\frak g$ is the $P$--submodule generated by $\Bbb
E_0\subset\ker(\square)$. In particular, $\Bbb E$ is a quotient of
$\oplus_{i\geq 0}(\otimes^i\frak p_+\otimes\bE_0)$. Likewise, for
$\al\in\Ga(F_1M)$, the section $L(\al)$ has values in a subbundle
associated to a quotient of $\oplus_{i\geq 0}(\otimes^i\frak
p_+\otimes\bF_1)$. Therefore, $\ka\wedge L(\al)$ has values in a
subbundle induced by a quotient of $\oplus_{i\geq 0}(\otimes^i\frak
p_+\otimes\bE_0\otimes\bF_1)$. 

{}From the point of view of $G_0$ there is no difference between
submodules and quotients. Hence if we form some semi--holonomic jet of
$\ka\wedge L(\al)$, it will be a section of a subbundle corresponding
to a representation, which, as a $G_0$--module, is contained in
$\oplus_{i\geq 0}(\otimes^i\frak p_+\otimes\bE_0\otimes\bF_1)$.
Since our assumptions imply that there is no nonzero
$G_0$--homomorphism from any such submodule to $\bF_2$, we can
complete the proof by showing that $D^2(\al)$ is obtained by applying
a strongly invariant operator to $\ka\wedge L(\al)$.

The latter fact has been proved in \cite{CD}, but for the sake of
completeness we give the simple argument: By definition, we have
$D(\al)=\pi_H(d^\nabla L(\al))$. Hence we may use $d^\nabla L(\al)$ as
a lift of $D(\al)$, so
$$
LD(\al)=d^\nabla L(\al)-Q\partial^*(\ka\wedge L(\al)).
$$
Applying $\pi_Hd^\nabla$ we conclude that
$$
D^2(\al)=\pi_H\o (\id-d^\nabla Q\partial^*)(\ka\wedge L(\al)).
$$
\end{proof}

\subsection{The Hasse graph and BGG diagrams in the complex case}\label{4.1}
To apply the vanishing Theorem systematically, we need to use certain
facts concerning the decomposition of $H_*(\frak p_+,\Bbb V)$ into
irreducible components. If $\frak g$ is complex, and $\Bbb V$ is a
complex irreducible representation, then $H_*(\frak p_+,\Bbb V)$ was
completely described as a $\frak g_0$--representation in
\cite{Kostant}. The answer is remarkably uniform and is expressed
using the Hasse graph and BGG diagrams. Let us briefly describe the
result. 

Let $\frak g$ be a complex semisimple Lie algebra, $\frak
h\subset\frak g$ a Cartan subalgebra and $\Delta\subset\frak h^*$ the
corresponding set of roots. Then the real subspace $\frak
h_0\subset\frak h$ on which all roots have real values is a real form
of $\frak h$ and the Killing form induces a positive definite inner
product on $\frak h_0$. Fix a choice $\De_+\subset\De$ of a positive
subsystem and let $\De_0$ be the corresponding set of simple roots.
For $\al\in\De$ let $\si_{\al}:\frak h_0\to\frak h_0$ denote the
reflection in the hyperplane $\ker(\al)$. The Weyl group $W$ of $\frak
g$ is the finite subgroup of the orthogonal group $O(\frak h_0)$
generated by these reflections. Then also the reflections $\si_{j}$
corresponding to simple roots $\al_j\in\De_0$ generate the group $W$.
For $w\in W$, the length $|w|$ is defined as the smallest positive
integer $n$ such that $w=\si_{j_1}\o \ldots\o\si_{j_n}$. Apart from
the evident action of $W$ on $\frak h_0^*$, which we write by
$\la\mapsto w(\la)$, there is also the \textit{affine action}. This is
defined by $w\cdot \la:=w(\la+\de)-\de$ for $\la\in\frak h_0^*$, where
$\de$ denotes half the sum of all positive roots.

The choices of $\frak h$ and $\De_+$ give rise to the standard Borel
subalgebra $\frak b\subset\frak g$, which is the sum of $\frak h$ and
all positive root spaces. The \textit{Hasse graph} for $\frak b$ (see
\cite{BGG}) is the directed graph with vertices the elements of $W$
and labeled arrows defined by $w\overset{\al}\longrightarrow w'$ if
and only if $|w'|=|w|+1$ and $\al\in\De_+$ is such that $w'=\si_\al
\circ w$.

Now if $\frak p\subset\frak g$ is obtained from some $|k|$--grading
then one can always choose $\frak h$ and $\De_+$ in such a way that
$\frak b\subset\frak p$. Moreover, $\frak p$ then is the standard
parabolic subalgebra corresponding to a subset $\Sigma\subset
\De_0$. Explicitly, a simple root lies in $\Si$ if and only if the
corresponding root space is contained in $\frak g_1$. The root spaces
of the other simple roots then lie in $\frak g_0$. We also split
$\De_+=\De_+(\frak g_0)\sqcup\De_+(\frak p_+)$ according to positive
root spaces being contained in the indicated subalgebras. 

Associating to every $w\in W$ the subset
$\Ph_w:=\{\al\in\Delta_+:w^{-1}(\al)\in-\Delta_+\}\subset\Delta_+$,
one obtains a bijection between $W$ the set of all subsets of
$\Delta_+$ having a certain property. The Hasse graph of $\frak p$ is
then defined as the subgraph of the Hasse graph of $\frak b$ consisting
of all vertices $w\in W^\frak p:=\{w\in W:\Ph_w\subset\Delta_+(\frak
p_+)\}$ and all edges connecting these vertices. It turns out that
only elements of $\Delta_+(\frak p_+)$ can occur as labels for the
remaining arrows. 
 
There is an alternative characterization of $W^\frak p$: We say that a
weight $\la\in\frak h_0^*$ is $\frak g$--dominant (respectively $\frak
p$--dominant) if $\langle\la,\al\rangle\geq 0$ for all $\al\in\De_0$
(respectively all $\al\in \De_0\setminus\Sigma$). Then $w\in W^\frak
p$ if and only if for one (or equivalently any) $\frak g$--dominant
weight $\la$ the weight $w\cdot\la$ is $\frak p$--dominant. Given a
$\frak g$--dominant weight $\la$, we define the \textit{BGG diagram}
associated to $\la$ to be the graph obtained from the Hasse graph of
$\frak p$ by replacing the edge $w$ by $w\cdot\la$. 

Following the conventions in \cite{BE}, will label representations by
the highest weight of the dual rather than the highest weight of the
given representation. Equivalently, $\bV=\bV_{\la}$ if $-\la$ is the
lowest weight of $\bV$. We will use this notation for representations
of both $G$ and $G_0$. With this convention, the vertices in the BGG
diagram associated to $\frak g$--dominant integral weight $\la,$ are
exactly the labels of the irreducible components (with respect to
$G_0$) of $H_*(\frak p_+,\Bbb V_{\la})$. The irreducible components of
$H_k(\frak p_+,\Bbb V_\la)$ are exactly the vertices corresponding to
elements $w\in W^\frak p$ such that $|w|=k$. The paper \cite{Kostant}
even gives an explicit description of a highest weight vector in each
component, which is often helpful in determining how the components
sit inside of $\La^*\frak p_+\otimes\Bbb V$.

There are efficient methods to work out the form of the labeled Hasse
diagrams for low gradings (for details see \cite{KS}).  We shall quote
the results in cases of interests below. Now we can formulate the
condition in the vanishing theorem \ref{3.3} in terms of weights as
follows.

\begin{lem*}
  Let $\bE_0$ be an irreducible component of $H_2(\frak p_+,\frak g)$
  and let $\bF_1=\bF_{\la_1}\subset H_k(\frak p_+,\Bbb V)$ and
  $\bF_2=\bF_{\la_2}\subset H_{k+2}(\frak p_+,\Bbb V)$ be irreducible
  components.
 
  If $\bF_2$ is isomorphic to a $G_0$-submodule of $\oplus_{i\geq
    0}\left(\otimes^i\frak p_+\otimes\bE_0\otimes \bF_1\right)$,
  then $\la_2-\la_1$ is a weight of $\oplus_{i\geq
    0}\left(\otimes^i\g_-\otimes\bE_0^*\right)$.
\end{lem*}
\begin{proof}
 First note that
  $\bF_2$ is  isomorphic  to a $G_0$-submodule of
$$
\left(\oplus_{i\geq 0}\otimes^i\frak p_+\right)\otimes\bE_0\otimes \bF_1,
$$
if and only 
if $\bF_2^*$ is  isomorphic  to a $G_0$-submodule of
$$
\left(\oplus_{i\geq 0}\left(\otimes^i\frak g_-\right)\otimes
\bE_0^*\otimes \bF_1^*\right).
$$
 
Now the result follows from the following fact which is well known in
representation theory of semisimple and reductive Lie algebras: Let
$\Bbb V$ be an irreducible representation of highest weight $\la$ and
$\Bbb W$ an arbitrary finite dimensional representation. Then the
highest weight of any irreducible component of $\Bbb V\otimes\Bbb W$
can be written as the sum of $\la$ and some weight of $\Bbb W$. 
\end{proof}

\subsection{Hasse graph and BGG diagrams in the real case}\label{4.2}
A $|k|$--grading on a real semisimple Lie algebra $\frak g$ induces a
$|k|$--grading on the complexification $\frak g^{\Bbb C}$. The subalgebras
$\frak p_+\subset\frak p\subset\frak g$ complexify to their counterparts
obtained from the complex $|k|$--grading. Using this, we can deduce
the decomposition of $H_k(\frak p_+,\Bbb V)$ from the complex case
discussed in \ref{4.1} above. 

Let us review some basic facts about representations of real Lie
algebras, see \cite{Oniscik,Silhan} for details. For a real Lie
algebra $\frak a$, a complex representation can be simply viewed as a
real representation $\Bbb V$ together with an $\frak a$--invariant
complex structure $J:\Bbb V\to\Bbb V$. In this case also $-J$ is an
$\frak a$--invariant almost complex structure and the resulting
representation of $\frak a$ is called the \textit{conjugate} of $\Bbb
V$ and denoted by $\bar{\Bbb V}$. A complex representation of $\frak a$
uniquely extends to a representation of the complexification $\frak
a^{\Bbb C}$, but on the level of $\frak a^{\Bbb C}$ the relation
between the representations $\Bbb V$ and $\bar{\Bbb V}$ is more involved
than on the level of $\frak a$.

If $\Bbb V$ is a real irreducible representation of $\frak a$, then
one can form the complexification $\Bbb V^\Bbb C$. If $\Bbb V$ does
not admit an $\frak g$--invariant complex structure, then $\Bbb V^{\Bbb
  C}$ is again irreducible. However if $\Bbb V$ does admit an
invariant complex structure, then $\Bbb V^{\Bbb C}\cong\Bbb
V\oplus\bar{\Bbb V}$.

Let us return to the question of decomposing $H_*(\frak p_+,\Bbb V)$
in the case of real $\frak g$. If $\Bbb V$ is a complex representation
of $\frak g$, then $\La^k\frak p_+\otimes\Bbb V\cong\La^k\frak
p_+^{\Bbb C}\otimes_{\Bbb C}\Bbb V$ and this is compatible with the
differentials. This easily implies that $H_*(\frak p_+,\Bbb V)$ is
simply the restriction to $\frak g_0\subset\frak g_0^{\Bbb C}$ of
$H_*(\frak p_+^{\Bbb C},\Bbb V)$. In particular, we obtain the same
decomposition into irreducibles as in the complex case.

On the other hand, let us assume that the representation $\Bbb V$ does
not admit a $\frak g$--invariant complex structure. Then $\Bbb V^{\Bbb
  C}$ is irreducible, and one easily shows that $H_*(\frak p_+^{\Bbb
  C},\Bbb V^{\Bbb C})$ is the complexification of $H_*(\frak p_+,\Bbb
V)$. For a $\frak p$--dominant weight $\mu$ let $\bar\mu$ be the
weight characterized by $\Bbb E_{\bar\mu}=\overline{\Bbb E_\mu}$. Now
for an irreducible component $\Bbb E_{\mu}\subset H_k(\frak p_+^{\Bbb
  C},\Bbb V^{\Bbb C})$ there are two possibilities. Either it is the
complexification of a real irreducible component in $H_*(\frak
p_+,\Bbb V)$. The other possibility is that also $\Bbb
E_{\bar\mu}\subset H_k(\frak p_+^{\Bbb C},\Bbb V^{\Bbb C})$ and there
is a complex irreducible component in $H_*(\frak p_+,\Bbb V)$ whose
complexification is $\Bbb E_{\mu}\oplus\Bbb E_{\bar\mu}$. It is shown
in \cite{Silhan} that the second possibility occurs if and only if
$\mu\neq\bar\mu$. 

Hence we obtain a complete description of the decomposition of
$H_*(\frak p_+,\Bbb V)$ in the real case. If $\Bbb V$ is complex and
has highest weight $\la$ as a representation of $\frak g^{\Bbb C}$
then the decomposition is given by the BGG diagram associated to
$\la$. An entry $\mu$ in this diagram corresponds to the restriction
of the complex representation $\Bbb E_\mu$ to $\frak g_0\subset\frak
g_0^{\Bbb C}$. 

If $\Bbb V$ is real (i.e.~does not admit an invariant complex
structure) then let $\la$ be the highest weight of the $\frak g^{\Bbb
  C}$ representation $\Bbb V^{\Bbb C}$. Then the decomposition of
$H_*(\frak p_+,\Bbb V)$ is obtained by identifying in the BGG diagram
associated to $\la$ each weight $\mu$ with $\bar\mu$. If $\mu=\bar\mu$
then the vertex corresponds to a real irreducible component, while
vertices obtained by identifying $\mu$ with $\bar\mu\neq\mu$
correspond to complex irreducible components. In particular, we
immediately obtain the following real version of lemma \ref{4.1}: 

\begin{lem*}
  Let $\bV$ be a real irreducible $\gog$-module.  Let $\bE_0$ be an
  irreducible component of $H_2(\frak p_+,\frak g)$ and let $\Bbb
  F_1\subset H_k(\frak p_+,\Bbb V)$ and $\Bbb F_2\subset H_{k+2}(\frak
  p_+,\Bbb V)$ be complex irreducible components such that $\la_i$ and
  $\bar\la_i$ are the labels of the irreducible components of their
  complexifications. 
  
  If $\bF_2$ is isomorphic to a $G_0$-submodule of $\oplus_{i\geq
    0}\left(\otimes^i\frak p_+\otimes\bE_0\otimes \bF_1\right)$,
  then at least one of $\la_2-\la_1$, $\bar\la_2-\la_1$,
  $\la_2-\bar\la_1$, and $\bar\la_2-\bar\la_1$ is a weight of
$\oplus_{i\geq 0}\left(\otimes^i\g_-^{\bC}\otimes(\bE_0^{\bC})^*\right)$.
\end{lem*}

\subsection{}\label{5.1}
We next study the Hasse and BGG diagrams in the case relevant for
quaternionic and Grassmannian geometries. Let us consider the algebra
$\frak g= \frak s\frak l(n+2,\bC)$ with the grading
$\g=\g_{-1}\oplus\g_0\oplus\g_1$ corresponding to the Dynkin diagram
(with $n+1$ nodes)
$$
\bxbdb{\al_1}{\al_2}{\al_3}{a_{n+1}},
$$
where $\al_i=e_i-e_{i+1}$ are the simple roots in the standard
notation for the $A_n$ series.  The other positive roots for $\frak g$
are then given by $\be^{ij}:=\al_i+\ldots+\al_j=e_i-e_{j+1}$ with
$1\leq i<j\leq n+1$ and we put $\be^{ii}=\al_i$.

The semisimple part of $\frak g_0$ is $\frak s\frak l(2,\bC)\oplus
\frak s\frak l(n,\bC)$. We get $\Delta_+(\frak g_0)=\{\be^{11}\}\cup
\{\be^{ij}|3\leq i \leq j\leq n+1\}$ with the two sets corresponding
to the two summands, and $\Delta_+(\frak p_+)=\{\be^{1j}|2\leq
j\leq n+1\}\cup\{\be^{2j}|2\leq j\leq n+1\}$.

The Hasse graph can be computed by the methods of the book \cite{BE},
or, together with labels over the arrows, by the methods of \cite{KS}.
It has a triangular shape, whose form is easily seen from the case
$n=4$ given below.
$$
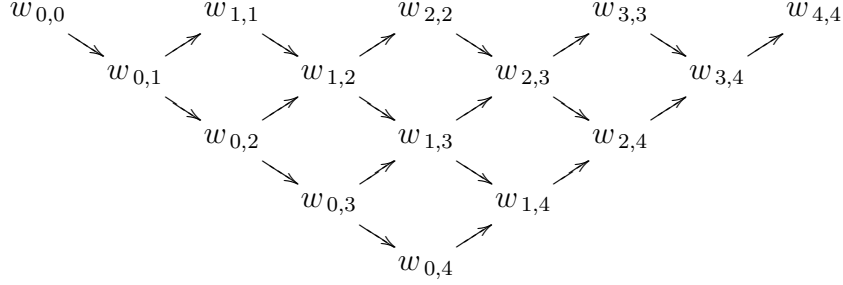
\xymatrix@R=3mm@C=3mm{%
{w_{\,0,0}}\ar[dr]& &{w_{\,1,1}}\ar[dr]& &{w_{\,2,2}}\ar[dr]
& &{w_{\,3,3}}\ar[dr]& &{w_{\,4,4}}
\\
&{w_{\,0,1}}\ar[ur]\ar[dr]&&{w_{\,1,2}}\ar[ur]\ar[dr]
&&{w_{\,2,3}}\ar[ur]\ar[dr]&&{w_{\,3,4}}\ar[ur]&
\\
&&{w_{\,0,2}}\ar[ur]\ar[dr]&&{w_{\,1,3}}\ar[ur]\ar[dr]&&{w_{\,2,4}}\ar[ur]&&
\\
&&&{w_{\,0,3}}\ar[ur]\ar[dr]&&{w_{\,1,4}}\ar[ur]&&&
\\
&&&&{w_{\,0,4}}\ar[ur]&&&&
}
$$

In general, the elements of length $k$ have the form $w_{i,j}$ with
$i+j=k$ and $0\leq i\leq j\leq n$. They can be computed explicitly,
but we do not need the result. Concerning the labels of the arrows,
the left edge of the diagram has the form
$$
w_{0,0}\overset{{\be^{22}}}\longrightarrow w_{0,1}
\overset{\be^{23}}\longrightarrow w_{0,2}
\overset{\be^{24}}\longrightarrow 
\ldots
\overset{\be^{2,n}}\longrightarrow w_{0,n-1}
\overset{\be^{2,n+1}}\longrightarrow w_{0,n},
$$
while the right edge of the diagram looks like
$$
w_{0,n}\overset{\be^{12}}\longrightarrow w_{1,n}
\overset{\be^{13}}\longrightarrow w_{2,n}
\overset{\be^{14}}\longrightarrow 
\ldots
\overset{\be^{1,n}}\longrightarrow w_{n-1,n}
\overset{\be^{1,n+1}}\longrightarrow w_{n,n}.
$$
In the rest of the diagram, parallel arrows have the same label.

For any representation $\Bbb V$ we will denote the corresponding
splitting of the Lie algebra homology groups as $H_k(\frak p_+,\Bbb
V)=\oplus H_{i,j}(\frak p_+,\Bbb V)$. In particular, we have
$H_2(\frak p_+,\frak g)=H_{0,2}(\frak p_+,\frak g)\oplus H_{1,1}(\frak
p_+,\frak g)$. 

\begin{prop*}
  Put $\Bbb E_0:=H_{1,1}(\frak p_+,\frak g)$. Then we have:

\noindent
(1) All weights of $\oplus_{i\geq
  0}\left(\otimes^i\g_-\otimes(\bE_0)^*\right)$ have the form
$$
-m_1\al_1-m_2\al_2+\textstyle{\sum_{i=3}^{n+1}} m_i\al_i
$$ 
for some integers $m_1,\dots,m_{n+1}$ such that $0<m_1<m_2$.

\noindent
(2) For all $i=0,\ldots,n-2;\,j=i,\ldots,n-2$ and any $\g$-dominant
  integral weight $\la$, we have
$$
w_{i,j+2}\cdot\la-w_{i,j}\cdot\la=
m_2\al_2+\ldots+m_{n+1}\al_{n+1}
$$
for some integers $m_2,\ldots,m_{n+1}$.

\noindent
(3) For all $j=2,\ldots,n ;\,i=0,\ldots,j-2$ and any
  $\g$-dominant integral weight $\la,$ we have 
$$
w_{i+2,j}\cdot\la-w_{i,j}\cdot\la=
-m(\al_1+\al_2)+m_3\al_3+\ldots+m_{n+1}\al_{n+1}
$$
for some integers $m,m_3\ldots,m_{n+1}$.
\end{prop*}
\begin{proof}
  (1) The standard recipes from \cite{BE} show that in terms of the
  fundamental weights $\la_i$ for $\frak g$ the highest weight of
  $\Bbb E_0^*$ is given by $-4\la_2+3\la_3+\la_{n+1}$. In particular, 
  the action of the semisimple part $\frak s\frak l(2,\Bbb
  C)\oplus\frak s\frak l(n,\Bbb C)$ of $\frak g_0$ is only via the
  second summand. Therefore, any weight of $\Bbb E_0^*$ is obtained by
  subtracting roots of the form $\be^{i,j}$ with $i\geq 3$ from the
  highest weight. In terms of simple roots, the highest weight of
  $\Bbb E_0^*$ reads as $-\al_1-2\al_2+\al_3+\dots+\al_{n+1}$, and hence
  any weight of $\Bbb E_0^*$ has the form
  $-\al_1-2\al_2+m_3\al_3+\dots+m_{n+1}\al_{n+1}$. On the other hand,
  the weights of $\frak g_-$ are exactly the elements of $-\De_+(\frak
  p_+)$. Since these are either of the form $-\al_1-\al_2-\dots-\al_j$
  or of the form $-\al_2-\dots-\al_j$, the claim follows.

\noindent
(2) If we have a labeled arrow $w\overset{\al}{\longrightarrow}w'$ in
    the Hasse diagram, then $w'=\si_\al(w)$, and hence for a $\frak
    g$--dominant weight $\la$ the difference $w'\cdot\la-w\cdot\la$ is
    an integer multiple of $\al$. Now from above we see that the Hasse
    diagram contains
    $w_{i,j}\overset{\be^{2,j+2}}{\longrightarrow}w_{i,j+1}
\overset{\be^{2,j+3}}{\longrightarrow}w_{i,j+2}$, and the claim
    follows since $\be^{2,\ell}=\al_2+\dots+\al_{\ell+1}$. 

\noindent
(3)  This is similar as in (2) taking into account that the Hasse
     diagram contains the part $w_{i,j}
\overset{\be^{1,i+2}}{\longrightarrow}w_{i+1,j}
\overset{\be^{1,i+3}}{\longrightarrow}w_{i+2,j}$ and that
     $\be^{1,\ell}=\al_1+\al_2+\dots+\al_\ell$. 
\end{proof}

\subsection{Grassmannian and quaternionic structures}\label{5.2}
There are two real forms of the grading considered in \ref{5.1} which
lead to well known geometric structures. Since we are dealing with a
$|1|$--grading here, an infinitesimal flag structure of type $(G,P)$
on $M$ (which is equivalent to a regular normal parabolic geometry,
see \ref{2.4}) is simply a first order $G_0$--structure, i.e.~a
reduction of the frame bundle of $M$ to the structure group $G_0$.

Putting $G=SL(n+2,\Bbb R)$, the subgroup $P$ turns out to be the
stabilizer of a plane and $G_0\cong S(GL(2,\Bbb R)\x GL(n,\Bbb
R))\subset GL(2n,\Bbb R)$. Hence these geometries exist on manifolds
of even dimension $2n$, and they are usually called \textit{almost
  Grassmannian structures}, see for example \cite[chapters 6 and
7]{AG} and \cite{Bailey-Eastwood} for the complex analog of these
geometries.  Essentially, such a structure is given by an isomorphism
from the tangent bundle $TM$ to the tensor product $E^*\otimes F$ for
two auxiliary bundles $E$ and $F$ of rank 2 and $n$, respectively.

The other choice of interest is $G=PSL(n+1,\Bbb H)$, so $\frak g$ is a
real form of $\frak s\frak l(2n+2,\Bbb C)$. Then $P$ turns out to be
the stabilizer of a quaternionic line in $\Bbb H^{n+1}$ and $G_0\cong
S(GL(1,\Bbb H)GL(n,\Bbb H))\subset GL(4n,\Bbb R)$. The resulting
geometry is called an \textit{almost quaternionic structure} on a
manifold $M$ of dimension $4n$, see \cite{Salamon}. It is given by a
rank 3subbundle $Q\subset L(TM,TM)$ which can be locally spanned by
$I$, $J$, and $IJ$ for two anti commuting almost complex structures
$I$ and $J$ on $M$. Lifting the structure to a two fold covering of
$G_0$ (which corresponds to replacing $G$ by $SL(n+1,\Bbb H)$ and
locally is uniquely possible), this is equivalent to a tensor product
decomposition of $TM\otimes\Bbb C$ into a factor of rank 2 and one of
rank $n$, which shows the similarity to the almost Grassmannian case.

Assume that we have a manifold $M$ equipped with one of these two
types of structures, $\Bbb V$ is an irreducible representation of $G$
and $VM$ is the corresponding tractor bundle. Then we have the bundles
$H_k(T^*M,VM)$ from \ref{2.4}, and they split according to the
decomposition of $H_k(\frak p_+,\Bbb V)$. One verifies directly that
the BGG diagram for $\frak g^{\Bbb C}$ as described in \ref{5.1} above
may never contain two conjugate weights, so by \ref{4.2} and \ref{5.1}
we always get $H_k(T^*M,VM)=\oplus H_{i,j}(T^*M,VM)$ with $i+j=k$ and
$1\leq i\leq j\leq n$ (respectively $2n$ in the Grassmannian case),
with the notation following \ref{5.1}. Restricting a BGG operator
$D$ to sections of one component $H_{i,j}(T^*M,VM)$ we obtain a
splitting $D=D^{1,0}+D^{0,1}$ with the two components having values in
sections of $H_{i+1,j}(T^*M,VM)$ and $H_{i,j+1}(T^*M,VM)$,
respectively. 

In particular, $H_2(T^*M,\Cal AM)=H_{0,2}(T^*M,\Cal AM)\oplus
H_{1,1}(T^*M,\Cal AM)$ and accordingly the harmonic curvature
decomposes into two parts. The part with values in $H_{0,2}(T^*M,\Cal
AM)$ in both cases can be determined as a specific component of the
torsion of an arbitrary linear connection on $TM$ which is compatible
with the $G_0$--structure. This component is independent of the choice
of the connection and it is the only part of the torsion that cannot
be eliminated by changing the connection. Hence it is exactly the
obstruction to torsion freeness in the sense of first order
structures, and its vanishing is also equivalent to torsion freeness
of the corresponding regular normal parabolic geometry. Torsion free
geometries are usually referred to as \textit{Grassmannian}
respectively \textit{quaternionic} structures.

\begin{thm*}
  Let $M$ be a smooth manifold of dimension $2n$ endowed with a
  Grassmannian structure or a quaternionic structure (which requires
  $n$ to be even). Let $\Bbb V$ be an irreducible representation of
  $G$ and let $VM$ be the corresponding tractor bundle. For $1\leq
  i\leq j\leq n$ put $\Cal H_{i,j}:=H_{i,j}(T^*M,VM)$. Then the BGG
  sequence associated to $VM$ contains the subcomplexes
\begin{align*}
  &\calH_{j,j}\overset{D^{0,1}}\longrightarrow \calH_{j,j+1}
  \overset{D^{0,1}}\longrightarrow \ldots 
  \overset{D^{0,1}}\longrightarrow \calH_{j,n}\qquad \text{ for
    }j=0,\dots,n-2\\
  &\calH_{0,j}\overset{D^{1,0}}\longrightarrow \calH_{1,j}
  \overset{D^{1,0}}\longrightarrow\ldots 
  \overset{D^{1,0}}\longrightarrow \calH_{j,j}\qquad \text{ for
    }j=2,\ldots,n
\end{align*}
\end{thm*}
\begin{proof}
Since the assumptions of Theorem \ref{3.3} are satisfied, and the BGG
diagrams have the same form as for the complexification, we can use
the weight condition from Proposition \ref{5.1}. For the compositions
$D^{0,1}\o D^{0,1}$, we see from part (2) of Proposition \ref{5.1}
that $\la_2-\la_1$ is a linear combination of $\al_2,\dots,\al_{n+1}$
only. For the composition $D^{1,0}\o D^{1,0}$ we see from part (3) of
Proposition \ref{5.1} that writing $\la_2-\la_1$ as a linear
combination of the $\al_i$, the roots $\al_1$ and $\al_2$ have the same
coefficient. Now the claim follows in both cases from part (1) of
Proposition \ref{5.1}.
\end{proof}

\subsection*{Remark} (1) In the quaternionic case, this result vastly
generalizes \cite{Salamon} and \cite{Baston}. Indeed, the complexes in
\cite{Salamon} are the $D^{0,1}$--complexes starting at $\Cal H_{0,0}$
in the special case that $\Bbb V$ is a symmetric power of the dual of
the standard representation. The paper \cite{Baston} contains the
$D^{1,0}$--complexes starting at $\Cal H_{0,n}$ for arbitrary $\Bbb
V$.

(2) The complexes constructed in the Theorem in general contain
    operators of arbitrarily high orders. For example in the
    $D^{0,1}$--complex starting at $\Cal H_{0,0}$ the orders look as
    follows: Suppose that $\Bbb V=\Bbb V_\la$ and
    $\la=a_1\la_1+\dots+a_{n+1}\la_{n+1}$ in terms of the fundamental
    weights. Then for each of the operators in the complex there is a
    unique $i$ such that the order is $a_i+1$. In particular, among
    these complexes the ones contained \cite{Salamon} are exactly
    those in which all operators are first order.

\subsection{}\label{5.3}
Next we study the case $\g=\frak s\frak l (n+2,\bC)$ with $n\geq 2$,
with the $|2|$--grading corresponding to the Dynkin diagram
$$
\xbdbx{\al_1}{\al_2}{}{\al_n}{\al_{n+1}}. 
$$
We continue to use the notation from \ref{5.1} for roots. The
semisimple part of $\frak g_0$ is $\frak s\frak l(n,\bC)$ and
$\De_+(\frak g_0)=\{\be^{ij}|2\leq i \leq j\leq n\}$. On the other
hand $\Delta_+(\frak p_+)$ contains all $\be^{1,j}$ and all
$\be^{i,n+1}$, and the root space of $\be^{1,n+1}$ coincides with
$\frak g_2$.

The general shape of the Hasse diagram can be seen from the example $n=3$:
$$
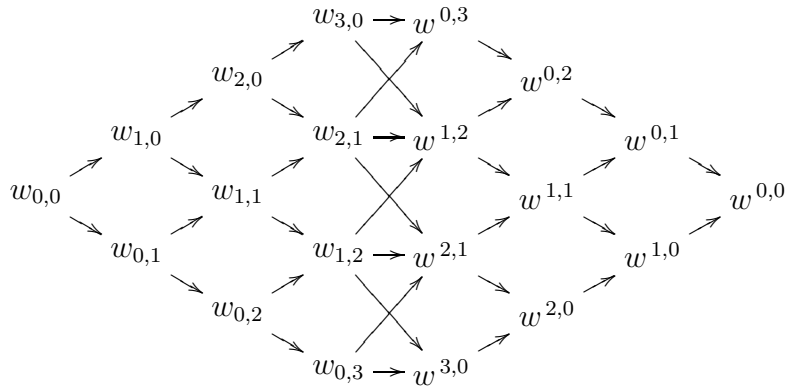
\xymatrix@R=2mm@C=4mm{%
&&&{w_{3,0}}
\ar[r]
\ar[ddr]
&{{w}^{\,0,3}}
\ar[dr]
&&&
\\
&&{w_{2,0}}
\ar[ur]
\ar[dr]
&&&{{w}^{0,2}}
\ar[dr]
&&
\\
&{w_{1,0}}
\ar[ur]
\ar[dr]
&& {w_{2,1}}
\ar[uur]
\ar[r]
\ar[ddr]
&{{w}^{\,1,2}}
\ar[ur]
\ar[dr]
&&{w}^{\,0,1}\ar[dr]
&
\\
w_{0,0}\ar[ur]\ar[dr]
&&w_{1,1}
\ar[ur]
\ar[dr]
&&&{{w}^{\,1,1}}
\ar[ur]\ar[dr]
&&
w^{\,0,0}
\\
&{w_{0,1}}
\ar[ur]
\ar[dr]
&& {w_{1,2}}
\ar[uur]
\ar[r]
\ar[ddr]
&{{w}^{\,2,1}}
\ar[ur]
\ar[dr]
&&{w}^{\,1,0}\ar[ur]
&
\\
&&{w_{0,2}}
\ar[ur]
\ar[dr]
&&&{{w}^{\,2,0}}
\ar[ur]
&&
\\
&&&{w_{0,3}}
\ar[r]
\ar[uur]
&{w}^{\,3,0}
\ar[ur]&&&
}
$$
Again the explicit form of the elements of $W^\frak p$ is not
important for our purposes. What we mainly need is that as in
\ref{5.1} parallel (or almost parallel) arrows have the same labels.
In particular, we have sequences of $n+1$ vertices each, which always
either go up or down. For the upward going sequences the labels over
the arrows are (in the right order) $\be^{1,1},
\be^{1,2},\dots,\be^{1,n}$, while for the downward going ones they are
$\be^{n+1,n+1},\be^{n,n+1},\dots,\be^{2,n+1}$. 

For any representation $\Bbb V$ of $\frak g$ and $k\leq n$ we
therefore obtain the splitting $H_k(\frak p_+,\Bbb V)=\oplus
H_{i,j}(\frak p_+,\Bbb V)$ with the sum over all $i,j\geq 0$ such that
$i+j=k$. For $k>n$ we obtain $H_k(\frak p_+,\Bbb V)=\oplus
H^{i,j}(\frak p_+,\Bbb V)$ with the sum over all $i,j\geq 0$ such that
$i+j=2n+1-k$. Similarly to Proposition \ref{5.1} one proves 

\begin{prop*}
Put $\bE_0:=H_{1,1}(\frak p_+,\frak g)$. Then we have:

\noindent
(1) All weights of $\oplus_{i\geq 0}\left(\otimes^i\g_-
  \otimes(\bE_0^{\bC})^*\right)$ have the form
$$
-m_1\al_1-m_{n+1}\al_{n+1}+\textstyle\sum_3^{n+1} m_i\al_i
$$
for integers $m_1,\dots,m_{n+1}$ such that $m_1,m_{n+1}>0$.

\noindent
(2) Let $\mu_1$ and $\mu_2$ be two weights which are contained in an
up going sequence in the BGG diagram of a $\frak g$--dominant integral
weight $\la$. Then $\mu_2-\mu_1$ can be written as a linear
combination of $\al_1,\dots,\al_n$. 

\noindent
(3) Let $\mu_1$ and $\mu_2$ be two weights which are contained in a
down going sequence in the BGG diagram of a $\frak g$--dominant
integral weight $\la$. Then $\mu_2-\mu_1$ can be written as a linear
combination of $\al_2,\dots,\al_{n+1}$.
\end{prop*}

\subsection{Lagrangean contact structures}\label{5.4}
There are various real forms of the situation discussed in \ref{5.3}
which are of interest in geometry. Putting $G:=SL(n+2,\Bbb R)$, one
obtains Lagrangean (or Legendrean) contact structures, see
\cite{Takeuchi} or \cite{C-tw}. Such a structure on a manifold $M$ of
dimension $2n+1$ is given by a codimension one subbundle $H\subset TM$
which defines a contact structure, and a fixed decomposition of
$H=E\oplus F$ as the direct sum of two Legendrean subbundles. This
means the the Lie bracket of two sections of $E$ (or two sections of
$F$) is a section of $H$.

Since we are dealing with a split real form here, all homology groups
(and hence also the corresponding vector bundles) split according to
the Hasse diagram discussed in \ref{5.3}. In particular, there are
three components in the harmonic curvature. The $(0,2)$-- and
$(2,0)$--parts are torsions which are the obstructions to
integrability of the subbundles $E,F\subset TM$. Vanishing of these
two components is equivalent to torsion freeness of the corresponding
parabolic geometry. Parallel to the proof of Theorem \ref{5.2},
Proposition \ref{5.3} leads to 

\begin{thm*}
Let $M$ a smooth manifold endowed with a torsion free Lagrangian
contact structure. Then the BGG sequence associated to any finite
dimensional irreducible representation of $\g$ splits according to the
Hasse diagram in \ref{5.3} and any upgoing or downgoing subsequence is
a complex.
\end{thm*}

\subsection{CR structures}\label{5.5}
The second class of interesting structures is obtained from
$G:=PSU(p+1,q+1)$ with $p\geq q$ and $p+q=n$. The resulting structures
are partially integrable almost CR structures of hypersurface type
which are non--degenerate of signature $(p,q)$, see \cite{CS}. The
analogy to Lagrangean contact structures can be seen by passing to the
complexified tangent bundle. In particular, starting with a complex
representation of $G$, the situation is completely parallel to the one
discussed in \ref{5.4}.

There is a difference however, in the case of real representations. If
$\Bbb V$ is a real representation and $\la$ is the highest weight of
$\Bbb V^{\Bbb C}$ then in the notation of \ref{5.3} one has
$\overline{(w_{i,j}\cdot\la)}=w_{j,i}\cdot\la$ and
$\overline{(w^{i,j}\cdot\la)}=w^{j,i}\cdot\la$, see \cite{Silhan}.
Hence the splitting of the real homologies $H_*(\frak p_+,\Bbb V)$ is
obtained from the Hasse diagram in \ref{5.3} by identifying the
$w_{i,j}$ with $w_{j,i}$ as well as $w^{i,j}$ with $w^{j,i}$.
Moreover, vertices with $i\neq j$ correspond to complex
subrepresentations in $H_*(\frak p_+,\Bbb V)$ while vertices with
$i=j$ correspond to real subrepresentations. The resulting picture for
$n=3$ looks as 
$$
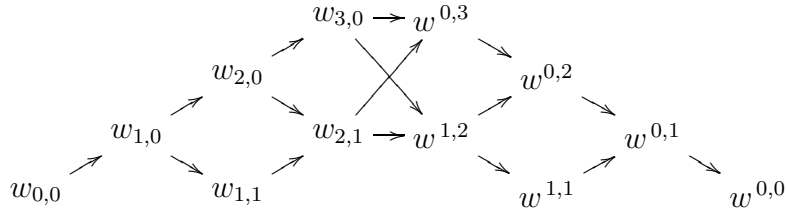
\xymatrix@R=2mm@C=4mm{%
&&&{w_{3,0}}
\ar[r]
\ar[ddr]
&{{w}^{\,0,3}}
\ar[dr]
&&&
\\
&&{w_{2,0}}
\ar[ur]
\ar[dr]
&&&{{w}^{0,2}}
\ar[dr]
&&
\\
&{w_{1,0}}
\ar[ur]
\ar[dr]
&& {w_{2,1}}
\ar[uur]
\ar[r]
%\ar[ddr]
&{{w}^{\,1,2}}
\ar[ur]
\ar[dr]
&&{w}^{\,0,1}\ar[dr]
&
\\
w_{0,0}\ar[ur]
&&w_{1,1}
\ar[ur]
&&&{{w}^{\,1,1}}
\ar[ur]
&&
w^{\,0,0}
}
$$

In particular, since the adjoint representation is real, there are
only two components in the harmonic curvature. There is just one
torsion, which is represented by the Nijenhuis tensor, and hence is
exactly the obstruction to integrability of the almost CR structure.
Integrable structures are usually referred to as CR structures.
 
\begin{thm*}
Let $M$ be smooth manifold endowed with a non--degenerate CR structure
of hypersurface type. 

\noindent (i) If $\bV$ is an complex irreducible representation of
$G$, then the associated BGG sequence splits according to the Hasse
diagram from \ref{5.3} and any upgoing or downgoing subsequence is a
complex.

\noindent (ii) If $\bV$ is a real irreducible representation of $G$,
then the associated BGG according to the diagram above, and any
upgoing or downgoing subsequence is a complex.
\end{thm*}

\begin{proof}
  The complex case is done as before. For the real case, we use Lemma
  \ref{4.2}. The differences $\la_2-\la_1$ and $\bar\la_2-\bar\la_1$
  can be handled as before. The differences $\bar\la_2-\la_1$ and
  $\la_2-\bar\la_1 $ also cannot be among the weights described in
  part (1) of Proposition \ref{5.3}, since in this difference either
  $\al_1$ or $\al_{n+1}$ must occur with a positive coefficient. 
\end{proof}

\subsection{}\label{5.6}
The last example we consider is $\g=\frak s\frak p(2k,\bC)$ for $k\geq
3$ with the $|2|$--grading described by the Dynkin diagram 
$$
\Cbxbdbb{\al_1}{\al_2}{\al_3}{\al_{k-1}}{\al_{k}}
$$ 
The positive roots are given by $\be^{i,j}=\al_i+\ldots+\al_j$ for
$1\leq i\leq j\leq k$ and $\ga^{i,j}=\be^{i,k-1}+\be^{jk}$ with $1\leq
i\leq k-1$ and $1\leq j\leq k$.  The semisimple part of $\g_0$ is
$\frak s\frak l(2,\bC)\oplus\frak s\frak p(2(k-2),\bC)$. We have 
$$
\De_+(\frak p_+)=\{\be^{i,j}:i=1,2;2\leq j\leq
k\}\cup\{\ga^{i,j}:i=1,2;1\leq j\leq k-1\}.
$$ 
The grading component $\frak g_2$ consists of the root spaces of
$\ga^{1,1}$, $\ga^{1,2}$, and $\ga^{2,2}$.  The general shape of the
Hasse diagram can be seen from the example $k=4$, which looks as
$$
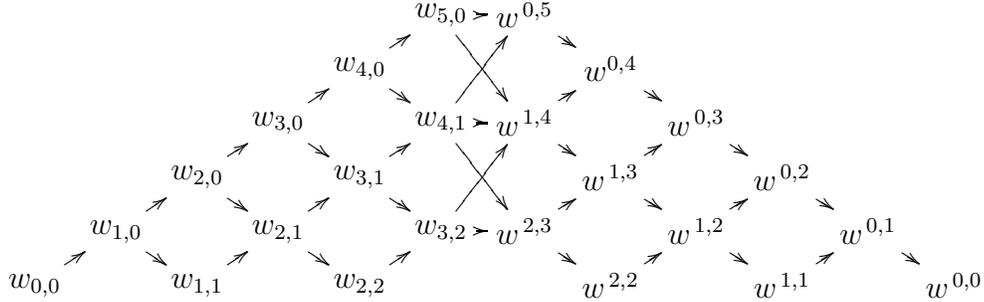
\xymatrix@R=1,4mm@C=1,4mm{%
&&&&&{w_{5,0}}
\ar[r]
\ar[ddr]
&{{w}^{\,0,5}}
\ar[dr]
&&&&&
\\
&&&&{w_{4,0}}
\ar[ur]
\ar[dr]
&&&{{w}^{0,4}}
\ar[dr]
&&&&
\\
&&&{w_{3,0}}
\ar[ur]
\ar[dr]
&& {w_{4,1}}
\ar[uur]
\ar[r]
\ar[ddr]
&{{w}^{\,1,4}}
\ar[ur]
\ar[dr]
&&{w}^{\,0,3}\ar[dr]
&&&
\\
&&w_{2,0}\ar[ur]\ar[dr]
&&w_{3,1}
\ar[ur]
\ar[dr]
&&&{{w}^{\,1,3}}
\ar[ur]\ar[dr]
&&
w^{\,0,2}\ar[dr]
&&
\\
&w_{1,0}\ar[ur]\ar[dr]
&&{w_{2,1}}
\ar[ur]
\ar[dr]
&& {w_{3,2}}
\ar[uur]
\ar[r]
 &{{w}^{\,2,3}}
\ar[ur]
\ar[dr]
&&{w}^{\,1,2}\ar[ur]\ar[dr]
&&w^{\,0,1}\ar[dr]&
\\
w_{0,0}\ar[ur]
&&w_{1,1}\ar[ur]
&&w_{2,2}\ar[ur]
&&&w^{\,2,2}\ar[ur]
&&w^{\,1,1}\ar[ur]
&&w^{\,0,0}
}
$$
 
For general $k$, the left edge of the diagram, including the labels
of the arrows, has the form 
$$
w_{0,0}\overset{{\be^{2,2}}}\longrightarrow 
\ldots
\overset{\be^{2,k-1}}\longrightarrow w_{k-2,0}
\overset{\ga^{2,2}}\longrightarrow w_{k-1,0}
\overset{\be^{2,k}}\longrightarrow
w_{k,0}\overset{\ga^{2,k-1}}\longrightarrow 
\ldots
\overset{\ga^{2,3}}\longrightarrow w_{2k-3,0}.
$$
The right edge of the diagram has the form 
$$
w^{0,2k-3}\overset{{\be^{1,2}}}\longrightarrow 
\ldots
\overset{\be^{1,k-1}}\longrightarrow w^{0,k-1}
\overset{\ga^{1,1}}\longrightarrow w^{0,k-2}
\overset{\be^{1,k}}\longrightarrow w^{0,k-3} 
\overset{\ga^{1,k-1}}\longrightarrow\ldots
\overset{\ga^{1,3}}\longrightarrow w^{0,0}.
$$
In the rest of the diagram, parallel (or almost parallel) arrows have
the same labels. As before, we will use the notation suggested by the
diagram for the irreducible components of the homology groups
$H_k(\frak p_+,\Bbb V)$.

\begin{prop*}
Put $\bE_0:=H_{1,1}(\frak p_+,\frak g)$. Then we have:

\noindent
(1) All weights of $\oplus_{i\geq 0}\left(\otimes^i(\g_-)
\otimes\bE_0^*\right)$ have the form 
$$
-m_1\al_1-m_2\al_2+\textstyle\sum_3^{k} m_i\al_i
$$
for integers $m_1,\ldots,m_{k}$ such that $0<m_1<m_2$.

\noindent
(2) Let $\mu_1$ and $\mu_2$ be two weights contained in an upgoing
 subsequence of the BGG diagram associated to some $\frak g$--dominant
 integral weight $\la$. Then the difference $\mu_2-\mu_1$ can be
 written as a linear combination of $\al_2,\dots,\al_k$. 

\noindent
(3) Let $\mu_1$ and $\mu_2$ be two weights contained in an downgoing
 subsequence of the BGG diagram associated to some $\frak g$--dominant
 integral weight $\la$. Then writing $\mu_2-\mu_1$ as a linear combination
 of $\al_1,\dots,\al_k$ the roots $\al_1$ and $\al_2$ have the same
 coefficient. 
\end{prop*}
\begin{proof}
The algorithms of \cite{BE} show that in terms of fundamental weights
the highest weight of $\Bbb E_0^*$ is $-5\la_2+4\la_3$, so in
particular the $\frak s\frak l(2,\Bbb C)$--factor in $\frak g_0$ acts
trivially on $\Bbb E_0^*$. In terms of simple roots we obtain
$-5\la_2+4\la_3=-\al_1-2\al_2+2(\al_3+\dots+\al_{k-1})+\al_k$. Now (1)
follows as in the proof of Lemma \ref{5.1}. For (2) we only have to
observe that in all roots which occur as labels in upgoing sequences,
the coefficient of $\al_1$ is trivial. Part (3) follows from the fact
that in all roots occurring in downgoing sequences the coefficients of
$\al_1$ and $\al_2$ are the same. 
\end{proof}

\subsection{Quaternionic and split--quaternionic contact
structures}\label{5.7} 
There are several real forms of the situation considered in \ref{5.6}
which are of interest in geometry. In all theses cases a, regular
normal parabolic geometry is equivalent to a certain codimension three
distribution $H\subset TM$ on a manifold $M$ of dimension
$4k-5$. Recall that given such a distribution and putting $Q:=TM/H$,
the Lie bracket of vector fields induces a tensorial map $H\x H\to
Q$. For each $x\in M$ this makes $H_x\oplus Q_x$ into a nilpotent
graded Lie algebra. The parabolic geometries then correspond to the
case that for each $x$ this is isomorphic to $\frak g_-$.

For $G=PSp(p+1,q+1)$ with $p\geq q$ and $p+q=k-1$, one obtains for
$\frak g_-$ the quaternionic Heisenberg algebra given by a
quaternionic Hermitian form of signature $(p,q)$. In particular, for
$q=0$ the resulting geometries are exactly the quaternionic contact
structures introduced by O.~Biquard, see \cite{Biq,Biq2}. The interest
in these structures comes from the fact that they occur as conformal
infinities of quaternionic K\"ahler manifolds. For the real form
$G=PSp(2k,\Bbb R)$, one obtains for $\frak g_-$ the (uniquely
determined) split quaternionic Heisenberg algebra. For $k=3$ and hence
$\dim(M)=7$ the two types of rank 4 distributions obtained in this way
are exactly the two generic types. 
 
Concerning the structure of the harmonic curvature, the case $k=3$ and
hence $\dim(M)=7$ is special. There are two independent harmonic
curvature components, one of which is a torsion and one of which is a
curvature. Hence torsion freeness is a nontrivial condition. It turns
out that torsion freeness is also equivalent to existence of a twistor
space. 

On the other hand, for $k>3$, the component $H_{2,0}(\frak p_+,\frak
g)$ consists of maps which are homogeneous of degree zero, and hence
vanishing of the corresponding component of the harmonic curvature is
a consequence of regularity. 

It turns out that in all cases, the BGG diagrams of the
complexification never contain conjugate weights. As for the other
geometries we obtain.
\begin{thm*}
Let $M$ a smooth manifold of dimension $4k-5,\,k\geq 3$, endowed with
a quaternionic contact structure or its split quaternionic
analog. Assume further that this structure is torsion free if $k=3$. 

Then for any irreducible representation $\Bbb V$ of $G$ the associated
BGG sequence splits according to the Hasse diagram from \ref{5.6} and
any upgoing or downgoing subsequence is a subcomplex. 
\end{thm*}

\section{Ellipticity}\label{6} 
In this section, we want to show that many of complexes obtained in
Theorem \ref{5.2} are elliptic in the quaternionic case.  To do this
we first analyze their symbol sequences in the Grassmannian case.

\subsection{Symbol sequences}\label{6.1}
As in \ref{5.2}, natural vector bundles on almost Grassmannian
manifolds are associated to representations of the group $P$ and in
the case of irreducible representations one has to deal with
$G_0=S(GL(2,\Bbb R)\x GL(n,\Bbb R))$. The standard representations
$\Bbb E$ and $\Bbb F$ of the two factors correspond to the bundles $E$
and $F$, and $T^*M\cong E\otimes F^*$. If we have two representations
$\Bbb V$ and $\Bbb W,$ then the symbol of an $r$th order differential
operator $D:\Ga(VM)\to\Ga(WM)$ between sections of the corresponding
bundles is a bundle map $S^rT^*M\otimes VM\to WM$. In the case of an
invariant differential operator, this bundle map is induced by a
$G_0$--equivariant map $\si:S^r(\Bbb E\otimes\Bbb F^*)\otimes\Bbb
V\to\Bbb W$. Determining all possible maps of this type is a sometimes
tedious but standard task in representation theory. For $X\in\Bbb
E\otimes\Bbb F^*$ we will write $\si_X:\Bbb V\to\Bbb W$ for the map
$v\mapsto \si(X\vee\dots\vee X\otimes v)$.

As a preliminary step, we have to analyze the representations $\Bbb
V^{k,\ell}:=S^k\Bbb E\otimes\La^\ell\Bbb F^*$. Then there is a unique
(up to scale) $G_0$--homomorphism $\si:(\Bbb E\otimes\Bbb
F^*)\otimes\Bbb V^{k,\ell}\to \Bbb V^{k+1,\ell+1}$ induced by the
symmetric product in the first, and the wedge product in the second
factor. Choosing a basis $\{e_1,e_2\}$ for $\Bbb E$, one may write any
element $X\in\Bbb E\otimes\Bbb F^*$ as $e_1\otimes
\al_1+e_2\otimes\al_2$ for elements $\al_1,\al_2\in\Bbb F^*$.
\begin{lem*}
For each $k\geq 0$ the symbol sequence 
$$
0\to\Bbb V^{k,0}\overset{\si_X}{\to}\Bbb V^{k+1,1}\to\dots
\overset{\si_X}{\to}\Bbb V^{k+n,n}\to 0
$$
is exact for $X=e_1\otimes \al_1+e_2\otimes\al_2$ provided that
$\al_1$ and $\al_2$ are linearly independent. 
\end{lem*}
\begin{proof}
  Let us assume throughout the proof that $\al_1$ and $\al_2$ are
  linearly independent. First consider the sequence
  $\Bbb V^{0,j-1}\overset{\si_X}{\to}\Bbb V^{1,j}\overset{\si_X}{\to}
  \Bbb V^{2,j+1}$.
  We claim that this sequence is exact for all $j=1,\dots,n-1$, the
  first map is injective for $j=1$ and the last map is surjective for
  $j=n-1$. 

Injectivity of the first map for $j=1$ is obvious. An arbitrary
element of $\Bbb V^{1,j}$ can be written as
$e_1\wedge\ph_1+e_2\wedge\ph_2$ for $\ph_1,\ph_2\in\La^j\Bbb
F^*$. Applying $\si_X$, we obtain
$$
e_1\vee e_1\otimes \al_1\wedge\ph_1+e_1\vee
e_2\otimes(\al_1\wedge\ph_2-\al_2\wedge\ph_1)+e_2\vee e_2\otimes
\al_2\wedge\ph_2.
$$ 
{}From this, surjectivity of the last map for $j=n-1$ follows
easily. Moreover, if this expression vanishes, then for $i=1,2$
vanishing of the coefficient of $e_i\vee e_i$ implies
$\al_i\wedge\ph_i=0$, so $\ph_i=\al_i\wedge\ps_i$. Vanishing of the
coefficient of $e_1\vee e_2$ leads to
$\al_1\wedge\al_2\wedge(\ps_2+\ps_1)=0$. Since $\al_1$ and $\al_2$ are
linearly independent, we obtain
$\ps_1+\ps_2=\al_1\wedge\rho_1+\al_2\wedge\rho_2$, and thus
$\ps_1-\al_1\wedge\rho_1=\ps_2-\al_2\wedge\rho_2=:\be$. 
By construction, $\al_i\wedge\be=\al_i\wedge\ps_i=\ph_i$ and hence
$e_1\otimes\ph_1+e_2\otimes\ph_2=\si_X(\be)$. 

Let us inductively assume that $k>1$ and we have proved that $\Bbb
V^{i-1,j-1}\overset{\si_X}{\to}\Bbb V^{i,j}\overset{\si_X}{\to} \Bbb
V^{i+1,j+1}$ is exact for all $i\leq k$ with the first map injective for
$j=1$ and the last map surjective for $j=n-1$. Consider the sequence 
$$
0\to S^{\ell-1}\Bbb E\otimes\La^2\Bbb E\to S^\ell\Bbb E\otimes\Bbb
E\to S^{\ell+1}\Bbb E\to 0,
$$
where the maps are given by symmetrization in the first $\ell$
respectively in all factors. Clearly the composition of these two maps
is trivial and looking at the dimensions one concludes that this is a
short exact sequence. The two maps are evidently compatible with
taking the symmetric product with some fixed element. Thus tensorizing
with appropriate exterior powers of $\Bbb F^*$ we obtain the following
commutative diagram with short exact columns, in which all horizontal
maps are either $\si_X$ or the tensor product of $\si_X$ with an
appropriate identity map:
$$
\xymatrix@R=5mm@C=5mm{%
\dots\to\Bbb V^{k-2,j-1}\ar[d]\ar[r] &\Bbb V^{k-1,j}\otimes\La^2\Bbb
 E\ar[d]\ar[r] & \Bbb V^{k,j+1}\otimes\La^2\Bbb E\to\cdots\ar[d]\\
\dots\to\Bbb V^{k-1,j-1}\otimes\Bbb E\ar[d]\ar[r]&\Bbb V^{k,j}\otimes\Bbb E\ar[d]\ar[r]&\Bbb V^{k+1,j+1}\otimes\Bbb E\to\cdots\ar[d]\\
\dots\to\Bbb V^{k,j-1}\ar[r]&\Bbb V^{k+1,j}\ar[r]&\Bbb V^{k+2,j+1}\to\cdots\\}
$$
By induction, the two top rows are exact, so exactness of the
bottom row (including the statements for $j=1$ and $j=n-1$) follows
from the nine--lemma of category theory.
\end{proof}

We are interested in the BGG sequences associated to the
representations $(S^k\Bbb V^*\otimes S^\ell\Bbb V)_0$, where $\Bbb
V=\Bbb R^{n+2}$ is the standard representation of $SL(n+2,\Bbb R)$ and
the subscript denotes the totally tracefree part. These are exactly
those representations of $G$ whose highest weight is a linear
combination of the first and last fundamental weights. We consider the
$D^{0,1}$--subcomplexes starting at $\Cal H_{0,0}$ (i.e.~the left edge
of the diagram) from Theorem \ref{5.2}. Let $\Bbb W^{k,\ell}_j$ be the
representation inducing $\Cal H_{0,j}$ for the given choice. It will
be convenient to put $\Bbb W^{k,\ell}_j=0$ for $j<0$ and $j>n$.

{}From the algorithms for determining Weyl orbits of weights in
\cite{BE} and the shape of the Hasse diagram, one
immediately concludes that $\Bbb W^{k,\ell}_j$ is the irreducible
component of highest weight in $S^k\Bbb E^*\otimes\Bbb W^{0,\ell}_j$
for all $j=0,\dots,n$. For $j<n$, one similarly concludes that $\Bbb
W^{k,\ell}_j$ is the irreducible component of highest weight in $\Bbb
W^{k,0}_j\otimes S^\ell\Bbb F$. Finally, one easily verifies directly
that $W^{0,0}_j=S^j\Bbb E\otimes\La^j\Bbb F^*\subset\La^j(\Bbb
E\otimes\Bbb F)$ for all $j$. Thus we conclude that 
$$
\Bbb W^{k,\ell}_j=(S^j\Bbb E\otimes S^k\Bbb E^*)_0\otimes (\La^j\Bbb
F^*\otimes S^{\ell}\Bbb F)_0
$$
for $j<n$. 

Since $\Bbb E$ has dimension two, the wedge product induces an
isomorphism $\Bbb E^*\cong \Bbb E\otimes\La^2\Bbb E^*$. Following the
usual conventions for conformal weights we indicate tensor product
with the $k$th power of the line bundle $\La^2\Bbb E^*$ by adding the
symbol $[k]$. Likewise, adding $[-k]$ indicates a tensor product with
the $k$th power of $\La^2\Bbb E$. The the above isomorphism reads as
$\Bbb E^*\cong \Bbb E[1]$. We also obtain an isomorphism $S^j\Bbb
E\otimes S^k\Bbb E^*\cong (S^j\Bbb E\otimes S^k\Bbb E)[k]$ under which
the tracefree part corresponds to $S^{j+k}\Bbb E[k]$. Finally, one
verifies directly that 
$$
\Bbb W^{k,\ell}_n=S^{k+n+\ell}\Bbb E[k]\otimes\La^n\Bbb F^*.
$$

The operators in our subcomplex are all of first order, except for the
last one, which is of order $\ell+1$. For $j<n-1$ there evidently is a
unique (up to scale) $G_0$--homomorphism $\Bbb E\otimes\Bbb
F^*\otimes\Bbb W^{k,\ell}_j\to\Bbb W^{k,\ell}_{j+1}$ which is induced
by taking the symmetric product in the $\Bbb E$ component and the
wedge product in the $\Bbb F^*$--component. In the last step, the
symbol should be a homomorphism 
$$
S^{\ell+1}(\Bbb E\otimes\Bbb F^*)\otimes S^{k+n-1}\Bbb E\otimes 
(\La^{n-1}\Bbb F^*\otimes S^{\ell}\Bbb F)_0\to 
S^{k+n+\ell}\Bbb E[k]\otimes\La^n\Bbb F^*. 
$$
Looking at the $\Bbb E$--components we see that such a homomorphism
has to factorize through $S^{\ell+1}\Bbb E\otimes S^{\ell+1}\Bbb
F^*\subset S^{\ell+1}(\Bbb E\otimes\Bbb F^*)$. But then there is again
a unique (up to scale) $G_0$--homomorphism. This is induced by the
symmetric product in the $\Bbb E$ component, while in the $\Bbb F^*$
component one has to take the unique contraction $S^{\ell+1}\Bbb
F^*\otimes S^{\ell}\Bbb F\to\Bbb F^*$ followed by the wedge product.
Hence we see that the symbols of the operators in the subcomplex are
all uniquely determined up to scale by their invariance properties.

\begin{thm*}
For all integers $k$ and $\ell$, the symbol sequence 
$$
0\to\Bbb W^{k,\ell}_0\overset{\si_X}{\to}\Bbb W^{k,\ell}_1\to\dots
\overset{\si_X}{\to}\Bbb W^{k,\ell}_n\to 0
$$
is exact for $X=e_1\otimes \al_1+e_2\otimes\al_2$ provided that
$\al_1$ and $\al_2$ are linearly independent.
\end{thm*}
\begin{proof}
  We proceed by induction on $\ell$. For $\ell=0$, we have $\Bbb
  W^{k,0}_j\cong \Bbb V^{k+j,j}[k]$, so the result follows directly
  from the Lemma.

Assuming that $\ell>0$, consider for $j=1,\dots,n-1$ the sequence   
$$
0\to (\La^{j-1}\Bbb F^*\otimes S^{\ell-1}\Bbb F)_0\to 
\La^j\Bbb F^*\otimes S^{\ell}\Bbb F
\to(\La^j\Bbb F^*\otimes S^{\ell}\Bbb F)_0\to 0,
$$
in which the first map is given by tensorizing with the identity and
then symmetrizing in the $\Bbb E$--part and alternating in the $\Bbb
F$--part, and the second map is projection to the tracefree
part. The dimensions of the representations can be easily computed
using Weyl's formula, and this shows that the sequence is short
exact. 

Tensorizing this exact sequence with $S^{k+j}\Bbb E[k]$, we obtain,
for each $j$, a short exact sequence
$$
0\to \Bbb W^{k+1,\ell-1}_{j-1}[-1] \to \Bbb W^{k,0}_j\otimes S^\ell\Bbb
F\to \Bbb W^{k,\ell}_j\to 0. 
$$
One easily verifies directly, that for $j<n-1$, we get a commutative diagram 
$$
\xymatrix@R=5mm@C=5mm{%
\Bbb E\otimes\Bbb F^*\otimes\Bbb W^{k+1,\ell-1}_{j-1}[-1]\ar[d]\ar[r] &
\Bbb W^{k+1,\ell-1}_j[-1]\ar[d]\\
\Bbb E\otimes\Bbb F^*\otimes\Bbb W^{k,0}_j\otimes S^\ell\Bbb F\ar[r] &
\Bbb W^{k,0}_{j+1}\otimes S^\ell\Bbb F}
$$
in which the vertical arrows come from the sequence above and the
horizontal arrows are tensor products of the symbol homomorphism $\si$
with appropriate identity maps. By exactness, these induce a
homomorphism $\Bbb E\otimes\Bbb F^*\otimes\Bbb W^{k,\ell}_j\to \Bbb
W^{k,\ell}_{j+1}$. From above, we know that this has to be a multiple
of $\si$, and this multiple has to be nonzero, since $\si\otimes\id$
maps onto $\Bbb W^{k,0}_{j+1}\otimes S^\ell\Bbb F$ by irreducibility
of $\Bbb W^{k,0}_{j+1}$.

Hence for $j=0,\dots,n-2$ we obtain a commutative diagram with short exact
columns in which the horizontal arrows are (nonzero multiples of)
$\si_X$ or the tensor product of $\si_X$ with an appropriate identity
map (recall that $\Bbb W^{r,s}_j=0$ for $j<0$):
$$
\xymatrix@R=5mm@C=5mm{%
\Bbb W^{k+1,\ell-1}_{j-2}[-1]\ar[r]\ar[d] & \Bbb
W^{k+1,\ell-1}_{j-1}[-1]\ar[r]\ar[d] & \Bbb W^{k+1,\ell-1}_j[-1]\ar[d]\\
\Bbb W^{k,0}_{j-1}\otimes\Bbb S^{\ell}F\ar[r]\ar[d] &
\Bbb W^{k,0}_j\otimes\Bbb S^{\ell}F \ar[r]\ar[d] & 
\Bbb W^{k,0}_{j+1}\otimes\Bbb S^{\ell}F\ar[d]\\ 
\Bbb W^{k,\ell}_{j-1} \ar[r] & \Bbb W^{k,\ell}_j \ar[r] & 
\Bbb W^{k,\ell}_{j+1}.}
$$
By induction the two top rows are exact, so exactness of our symbol
sequence at $\Bbb W^{k,\ell}_j$ for $j=0,\dots,n-2$ follows from the
nine lemma. 

For the last part, we first observe that $\Bbb
W^{k+1,\ell-1}_{n}[-1]\cong \Bbb W^{k,\ell}_{n}$. On the other hand,
for $j=n$, the above short exact sequences degenerate to isomorphisms
$(\La^{n-1}\Bbb F^*\otimes S^{\ell-1}\Bbb F)_0\cong\La^n\Bbb
F^*\otimes S^\ell\Bbb F$ respectively $\Bbb
W^{k+1,\ell-1}_{n-1}[-1]\cong\Bbb W^{k,0}_n\otimes S^\ell\Bbb F$.

Hence we obtain the following commutative diagram in which the first
two columns are short exact, and the horizontal arrows are nonzero
multiples of $\si_X$ or a tensor product of $\si_X$ with an
appropriate identity map:
$$
\xymatrix@R=5mm@C=4mm{%
\Bbb W^{k+1,\ell-1}_{n-3}[-1]\ar[r]\ar[d] & 
\Bbb W^{k+1,\ell-1}_{n-2}[-1] \ar[r]\ar[d] &
\Bbb W^{k+1,\ell-1}_{n-1}[-1]\ar[r]\ar[d]_\cong &
\Bbb W^{k+1,\ell-1}_{n}[-1]\to 0\\
\Bbb W^{k,0}_{n-2}\otimes S^{\ell}\Bbb F\ar[r]\ar[d] &
\Bbb W^{k,0}_{n-1}\otimes S^{\ell}\Bbb F\ar[r]\ar[d] &
\Bbb W^{k,0}_{n}\otimes S^{\ell}\Bbb F\ar[r] & 0\phantom{\to 0}\\
\Bbb W^{k,\ell}_{n-2}\ar[r] & \Bbb W^{k,\ell}_{n-1} && }
$$
By induction, the two top rows are exact. We can define a map
${\Bbb W}^{k,\ell}_{n-1}\to {\Bbb W}^{k+1,\ell-1}_{n}[-1]\cong {\Bbb
  W}^{k,\ell}_{n}$ as follows: Choose a preimage in ${\Bbb
  W}^{k,0}_{n-1}\otimes S^{\ell}\Bbb F$, map it to ${\Bbb
  W}^{k,0}_{n}\otimes S^{\ell}\Bbb F$ go up via the isomorphism, and
map to ${\Bbb W}^{k+1,\ell-1}_{n}[-1]$. Diagram chasing shows that
this is well defined and inserting it as the last map in the sequence
we get exactness at $\Bbb W^{k,\ell}_{n-1}$ and $\Bbb W^{k,\ell}_n$.
To complete the proof, it thus remains to show that this map is a
nonzero multiple of $\si_X$.

The inclusion of the tracefree part into $\La^{n-1}\Bbb F^*\otimes
S^\ell\Bbb F$ induces a $G_0$--homomorphism $\Bbb
W^{k,\ell}_{n-1}\to\Bbb W^{k,0}_{n-1}\otimes S^\ell\Bbb
F$. Tensorizing with the identity on $\Bbb E\otimes\Bbb F^*$ and
composing, we get a homomorphism
$$
\Bbb E\otimes\Bbb F^*\otimes\Bbb W^{k,\ell}_{n-1}\to 
\Bbb E\otimes\Bbb F^*\otimes\Bbb W^{k,0}_{n-1}\otimes S^\ell\Bbb
F\overset{\si\otimes\id}{\longrightarrow}\Bbb W^{k,0}_n\otimes
S^\ell\Bbb F.
$$  
{}From above we know that the target of this homomorphism is isomorphic
to $\Bbb W^{k+1,\ell-1}_{n-1}[-1]$ and hence in particular
irreducible. A moment of thought shows that the composition is nonzero
and hence surjective by irreducibility. Hence we have obtained a
surjective $G_0$--homomorphism $\Bbb E\otimes\Bbb F^*\otimes\Bbb
W^{k,\ell}_{n-1}\to\Bbb W^{k+1,\ell-1}_{n-1}[-1]$. Tensorize this with
the identity on $S^\ell(\Bbb E\otimes\Bbb F^*)$ and consider the
composition 
$$
S^\ell(\Bbb E\otimes\Bbb F^*)\otimes \Bbb E\otimes\Bbb F^*\otimes\Bbb
W^{k,\ell}_{n-1}\to S^\ell(\Bbb E\otimes\Bbb F^*)\otimes \Bbb
W^{k+1,\ell-1}_{n-1}[-1]\overset{\si}{\longrightarrow} \Bbb W^{k,\ell}_n. 
$$
This is surjective by induction, and looking at the explicit form
of the representation $\Bbb W^{k,\ell}_n$ one immediately sees that it
has to factor through $S^{\ell+1}\Bbb E\otimes S^{\ell+1}\Bbb
F^*\subset S^\ell(\Bbb E\otimes\Bbb F^*)\otimes \Bbb E\otimes\Bbb
F^*$. Therefore, it restricts to a nonzero multiple of the symbol map
on that part, and inserting copies of $X$, the claim follows.
\end{proof}

\subsection{Dual sequences}\label{6.10}
By duality, we can prove ellipticity for the $D^{1,0}$--subcomplexes
starting at $\Cal H_{0,n}$ in the BGG sequences considered in
\ref{6.1}. Let us denote by $\tilde{\Bbb W}^{k,\ell}_j$ the
$G_0$--representation corresponding to the bundle $\Cal H_{j,n}$ in
the BGG sequence associated to $(S^k\Bbb V^*\otimes S^\ell\Bbb V)_0$.
The crucial point here is that $\tilde{\Bbb W}^{k,\ell}_j\cong (\Bbb
W^{\ell,k}_{n-j})^*\otimes\La^{2n}\frak g_{-}^*$. 

This isomorphism comes from the bilinear map 
$$
\La^{n+j}\frak g_-^*\otimes (S^k\Bbb V^*\otimes S^\ell\Bbb V)_0\x
\La^{n-j}\frak g_-^*\otimes(S^k\Bbb V\otimes S^\ell\Bbb V^*)_0\to
\La^{2n}\frak g_-^*
$$
given by the wedge product and the paring between dual
representations. Note that $\dim(\frak g_-)=2n$, so $\La^{2n}\frak
g_-^*$ is one--dimensional.

In particular, this implies that dualizing and tensorizing with the
identity of this one--dimensional representation induces an
isomorphism 
$$
L(\tilde{\Bbb W}^{k,\ell}_j,\tilde{\Bbb W}^{k,\ell}_{j+1})\cong 
L(\Bbb W^{\ell,k}_{n-j-1},\Bbb W^{\ell,k}_{n-j}).
$$
In this case, the first operator of the sequence is of order $k+1$,
while all others are of first order. Now the symbol of an $r$th order
natural differential operator can equivalently be interpreted as a
$G_0$--equivariant map from $S^r(\Bbb E\otimes\Bbb F^*)$ to the
module of linear maps between the representations inducing the
bundles. 

Thus, the results of \ref{6.1} immediately imply that the symbols in
our sequence are uniquely determined up to scale by
$G_0$--equivariancy. Moreover, for any element $X\in \Bbb E\otimes\Bbb
F^*$, and each $j$, the symbol map $\si_X:\tilde{\Bbb W}^{k,\ell}_j\to
\tilde{\Bbb W}^{k,\ell}_{j+1}$ is the dual of the symbol map 
$\si_X:\Bbb W^{\ell,k}_{n-j-1}\to\Bbb W^{\ell,k}_{n-j}$. Since the
dual of an exact sequence is exact, we obtain

\begin{thm*}
The symbol sequence 
$$
0\to\tilde{\Bbb W}^{k,\ell}_0\overset{\si_X}{\to}\tilde{\Bbb
  W}^{k,\ell}_1\to\dots \overset{\si_X}{\to}\tilde{\Bbb W}^{k,\ell}_n\to 0
$$
is exact for $X=e_1\otimes \al_1+e_2\otimes\al_2$ if
$\al_1$ and $\al_2$ are linearly independent.
\end{thm*}

\subsection{Elliptic complexes for quaternionic structures} 
The results on symbol sequences in the Grassmannian case derived above
have immediate consequences for quaternionic structures:

 \begin{thm*}
Let $M$ be a quaternionic manifold of dimension $4n\geq 8$. 
Let $\Bbb V=\Bbb H^{n+1}$ be the standard representation of
$SL(n+1,\Bbb H)$. 

Then for all integers $k,\ell\geq 0$, the subcomplexes 
\begin{gather*}
0\to \Ga(\Cal H_{0,0}M)\overset{D^{0,1}}\to\Ga(\Cal H_{0,1}M)
\overset{D^{0,1}}\to\dots\overset{D^{0,1}}\to\Ga(\Cal H_{0,n}M)\to 0\\
0\to \Ga(\Cal H_{0,n}M)\overset{D^{1,0}}\to\Ga(\Cal H_{1,n}M)
\overset{D^{1,0}}\to\dots\overset{D^{1,0}}\to\Ga(\Cal H_{n,n}M)\to 0
\end{gather*}
of the BGG sequence associated to the representation $S^k\Bbb
V^*\otimes S^\ell\Bbb V$ are elliptic. In particular, this applies to
the deformation complex for quaternionic structures, see
\cite{deformations}, which is the $D^{0,1}$--complex for $k=\ell=1$. 
\end{thm*}
  
\begin{proof}
  The symbol sequences have the same complexifications as the ones in
  Theorems \ref{6.1} and \ref{6.10}. Since a nonzero quaternionic
  linear map $\Bbb H\to\Bbb H^n$ always has complex rank two, the
  condition for exactness of the symbol sequence in these Theorems is
  always satisfied.
\end{proof}
 
The $D^{0,1}$--complexes for $\ell=0$ are all the elliptic complexes
found in \cite{Salamon}, except the one for $r=-1$ in \cite[Theorem
5.5]{Salamon}, which belongs to singular infinitesimal character.

\end{document}